\documentclass[letterpaper, 10pt, journal, onecolumn, final]{IEEEtran}
\usepackage{subfig}
\usepackage{cite}
\usepackage{textcomp}

\usepackage{algorithm}
\usepackage{algpseudocode}
\usepackage{algorithmicx}
\usepackage{graphicx}

\usepackage{amssymb,amsmath,amsfonts,cases,mathtools}
\usepackage{microtype}
\usepackage{url}
\usepackage{hyperref}
\usepackage{xifthen}
\usepackage{xpatch}
\usepackage{amsthm}
\usepackage[dvipsnames]{xcolor}
\usepackage{tikz}
\usepackage{lipsum,booktabs}
\usepackage{pifont}
\usepackage{xcolor}
\usepackage{soul}
\usepackage{balance}

\newcommand{\R}{\mathbb{R}}

\newcommand{\N}{\mathbb{N}}

\newcommand{\T}{^\top}

\newcommand{\until}[1]{\{1,\ldots,#1\}}

\DeclareMathOperator{\col}{col}
\DeclareMathOperator{\blkdiag}{blkdiag}

\newtheorem{theorem}{Theorem}

\newtheorem{lemma}{Lemma}
\newtheorem{assumption}{Assumption}
\newtheorem{remark}{Remark}

\newcommand\oprocendsymbol{\hbox{$\blacksquare$}}
\newcommand\oprocend{\relax\ifmmode\else\unskip\hfill\fi\oprocendsymbol}

\newcommand{\cB}{\mathcal B}

\newcommand{\cE}{\mathcal E}

\newcommand{\cG}{\mathcal G}

\newcommand{\cK}{\mathcal K}

\newcommand{\cM}{\mathcal M}
\newcommand{\cN}{\mathcal N}

\newcommand{\cS}{\mathcal S}

\newcommand{\cX}{\mathcal X}

\graphicspath{{figs/}}

\newcommand{\norm}[1]{\left\lVert#1 \right\rVert}

\newcommand{\rz}{\mathrm{z}}

\newcommand{\rx}{\mathrm{x}}
\newcommand{\rw}{\mathrm{w}}

\newcommand{\uu}{G}
\newcommand{\z}{\rz}

\newcommand{\x}{\rx}
\newcommand{\w}{\rw}

\newcommand{\h}{h}
\newcommand{\bb}{b}
\newcommand{\g}{g}
\newcommand{\cc}{c}
\renewcommand{\aa}{a}

\newcommand{\ud}{^}
\newcommand{\du}{_}

\newcommand{\proxy}{\quad \text{proxy for} \quad}

\newcommand{\iter}{k}
\newcommand{\iterp}{{\iter+1}}

\newcommand{\initer}{\tau}

\newcommand{\tny}[1]{\text{\tiny #1}}
\newcommand{\av}{\du{\tny{AV}}}

\newcommand{\n}{n}

\newcommand{\zt}{\z\ud\iter}

\newcommand{\tx}{\tilde{\x}}
\newcommand{\txt}{\tx\ud\iter}
\newcommand{\txtp}{\tx\ud{\iterp}}

\newcommand{\f}{f}

\newcommand{\phii}{\phi\du{i}}
\newcommand{\phij}{\phi\du{j}}

\newcommand{\xstari}{\xstar\du{i}}

\newcommand{\1}{\mathbf{1}}

\newcommand{\zit}{\zt\du{i}}

\newcommand{\chit}{\chi\ud\iter}

\newcommand{\step}{\gamma}

\newcommand{\nc}{n_\chi}

\newcommand{\fs}{\f_{\sigma}}

\newcommand{\inputnn}{u}

\newcommand{\uxi}{\inputnn\du{x,i}} 
\newcommand{\usi}{\inputnn\du{\sigma,i}}

\newcommand{\set}{\until{N}}

\newcommand{\dr}{\mathrm{d}} %
\newcommand{\dt}{\dr\ud{\iter}}

\newcommand{\dta}{\dr\ud{\initer}}

\newcommand{\dxt}{\dr\ud{\iter}\du{x}} 
\newcommand{\dxit}{\dt\du{x,i}}

\newcommand{\dxita}{\dta\du{x,i}}

\newcommand{\dst}{\dr\ud{\iter}\du{\sigma}} 
\newcommand{\dsit}{\dt\du{\sigma,i}}

\newcommand{\dsita}{\dta\du{\sigma,i}}

\newcommand{\hfi}{\hat{\f}\du{i}}
\newcommand{\hf}{\hat{\f}}

\allowdisplaybreaks

\newcommand{\adjSymbol}{A}
\newcommand{\adj}{\mathcal{\adjSymbol}}
\newcommand{\adjkron}{\adjSymbol}
\newcommand{\adjentry}{a}

\newcommand{\adjminusI}{M}

\newcommand{\xii}{\x\du{i}}
\newcommand{\xit}{\x\ud\iter\du{i}}
\newcommand{\xjt}{\x\ud\iter\du{j}}
\newcommand{\xitp}{\x\ud{\iterp}\du{i}}

\newcommand{\xt}{\x\ud\iter}
\newcommand{\xtp}{\x\ud{\iterp}}
\newcommand{\xstar}{x\ud{\star}}

\newcommand{\txtav}{\txt\av}
\newcommand{\txav}{\tx\av}
\newcommand{\txtpav}{\txtp\av}

\newcommand{\stracker}{\w}

\newcommand{\si}{\stracker\du{i}}

\newcommand{\sit}{\stracker\ud\iter\du{i}}
\newcommand{\sjt}{\stracker\ud\iter\du{j}}
\newcommand{\sitp}{\stracker\ud\iterp\du{i}}

\newcommand{\ytracker}{\z}

\newcommand{\yit}{\ytracker\ud\iter\du{i}}
\newcommand{\yjt}{\ytracker\ud\iter\du{j}}
\newcommand{\yitp}{\ytracker\ud\iterp\du{i}}

\newcommand{\wit}{\stracker\ud\iter\du{i}}

\newcommand{\wt}{\stracker\ud\iter}
\newcommand{\wtp}{\stracker\ud\iterp}

\newcommand{\wm}{\bar{\stracker}}
\newcommand{\wmt}{\wm\ud\iter}
\newcommand{\wmtp}{\wm\ud\iterp}
\newcommand{\winitm}{\wm\ud 0}

\newcommand{\wort}{\stracker\du\perp}
\newcommand{\wortt}{\wort\ud\iter}
\newcommand{\worttp}{\wort\ud\iterp}

\newcommand{\twort}{\tilde{\stracker}\du{\perp}}
\newcommand{\twortt}{\twort\ud\iter}
\newcommand{\tworttp}{\twort\ud\iterp}

\newcommand{\qit}{\ytracker\ud\iter\du{i}}

\newcommand{\qt}{\ytracker\ud\iter}
\newcommand{\qtp}{\ytracker\ud\iterp}

\newcommand{\qm}{\bar{\ytracker}}
\newcommand{\qmt}{\qm\ud\iter}
\newcommand{\qmtp}{\qm\ud\iterp}
\newcommand{\qinitm}{\qm\ud 0}

\newcommand{\qort}{\ytracker\du{\perp}}
\newcommand{\qortt}{\qort\ud\iter}
\newcommand{\qorttp}{\qort\ud\iterp}

\newcommand{\tqort}{\tilde{\ytracker}\du{\perp}}
\newcommand{\tqortt}{\tqort\ud\iter}
\newcommand{\tqorttp}{\tqort\ud\iterp}

\newcommand{\hu}{\hat{\uu}}

\newcommand{\huxt}{\hu\du{1}}
\newcommand{\hust}{\hu\du{2}}

\newcommand{\hud}{\hat{\uu}\ud{d}}

\newcommand{\Lfs}{L}
\newcommand{\Lfx}{L\du{x}}
\newcommand{\Lfsigma}{L\du{\sigma}}
\newcommand{\Lfth}{L\du{\theta}}

\newcommand{\mufx}{\mu\du{\f}}

\newcommand{\Lhfx}{\hat{L}\du{x}}
\newcommand{\Lhfsigma}{\hat{L}\du{\sigma}}

\newcommand{\Lphi}{L\du{\phi}}

\newcommand{\LF}{L\du{\ell}}

\newcommand{\betath}{\beta_{\theta}}
\newcommand{\betax}{\beta_{x}}

\newcommand{\loss}{\ell}
\newcommand{\lossi}{\ell\du{i}}
\newcommand{\multh}{\mu\du{\loss}}
\newcommand{\Llth}{L\du{\ell}}

\newcommand{\lossav}{\tilde{\loss}_{\Tp}}

\newcommand{\nablalosst}{\nabla_3 \tilde{\loss}}
\newcommand{\nablalossav}{\nabla_2 \lossav}

\newcommand{\thh}{\theta}
\newcommand{\thi}{\theta\du{i}}
\newcommand{\thj}{\theta\du{j}}
\newcommand{\thit}{\thi\ud\iter}
\newcommand{\thjt}{\thj\ud\iter} 
\newcommand{\thitp}{\thi\ud\iterp}

\newcommand{\thistar}{\thi\ud\star}

\newcommand{\tth}{\tilde{\theta}}
\newcommand{\tht}{\theta\ud\iter} 
\newcommand{\thtp}{\theta\ud\iterp}
\newcommand{\thstar}{\theta\ud\star}
\newcommand{\ttht}{\tilde{\theta}\ud\iter}
\newcommand{\tthtp}{\tilde{\theta}\ud\iterp}

\newcommand{\tthtav}{\ttht\av}
\newcommand{\tthav}{\tth\av}
\newcommand{\tthtpav}{\tthtp\av}

\newcommand{\thSet}{\Theta}

\newcommand{\Ball}{\cB}
\newcommand{\thBall}{\Ball_{\bar{r}\du{i}}(\thistar)}

\newcommand{\dimtheta}{m} 
\newcommand{\dimthetai}{m\du{i}}

\newcommand{\Tp}{\cK}

\newcommand{\Rmat}{R}

\newcommand{\Vslow}{W}

\newcommand{\Vslowth}{\Vslow_{\thh}}
\newcommand{\Vslowtx}{\Vslow_{\x}}
\newcommand{\Vfast}{U}
\newcommand{\Vconv}{W}
\newcommand{\Vnom}{V}
\newcommand{\dVslow}{\Delta \Vslow}
\newcommand{\dVfast}{\Delta \Vfast}

\newcommand{\dVnom}{\Delta \Vnom}

\newcommand{\gSlow}{\g\du{\textsc{r}}}
\newcommand{\gSlowx}{\g\du{\textsc{r,2}}}
\newcommand{\gSlowth}{\g\du{\textsc{r,1}}}
\newcommand{\gFast}{\g\du{\textsc{b}}}

\newcommand{\ath}{\alpha}

\newcommand{\Pmatslow}{P\du{\Vslow}}
\newcommand{\Pmatfast}{P\du{\Vfast}}
\newcommand{\Plyapfast}{\Pi}

\newcommand{\Qmatfastentry}{q}

\newcommand{\RARdiag}{\tilde{\adjkron}}
\newcommand{\Kmat}{\cM(\Plyapfast)}

\newcommand{\domVpert}{D}

\newcommand{\stateF}{\zeta}
\newcommand{\stateFt}{\stateF\ud{\iter}}
\newcommand{\stateFtp}{\stateF\ud{\iterp}}

\newcommand{\stateS}{\chi}
\newcommand{\stateSt}{\stateS\ud{\iter}}
\newcommand{\stateStp}{\stateS\ud{\iterp}}

\newcommand{\stateSav}{\chi\av}

\newif\ifkappa
\kappatrue

\newcommand{\citeConverse}{\cite[Th.~2.2.1]{bai1988averaging}}

\newcommand{\state}{\xi}

\newcommand{\statet}{\state\ud\iter} 
\newcommand{\statetp}{\state\ud\iterp}

\newcommand{\dyn}{F_{\step}} 
\newcommand{\dyncomp}[1]{F_{\step, #1}} 
\newcommand{\dyno}{F_{\step}^{\tny{or}}} 
\newcommand{\dynocomp}[1]{F_{\step,#1}^{\tny{or}}} 

\newcommand{\nstate}{n_{\state}} 

\newcommand{\pert}{\mathrm{p}}

\newcommand{\perto}{\pert\du{1}}
\newcommand{\pertt}{\pert\du{2}}
\newcommand{\pertzi}{\pert\du{\loss,i}}
\newcommand{\pertoi}{\pert\du{1,i}}
\newcommand{\pertti}{\pert\du{2,i}}

\newcommand{\pertz}{\pert\du{\loss}}

\newcommand{\bound}{B}

\newcommand{\pp}{p}
\newcommand{\aat}{\tilde{\aa}}

\newcommand{\quadd}{P}
\newcommand{\quadi}{\quadd_i}
\newcommand{\lin}{v}
\newcommand{\lini}{\lin_i}
\newcommand{\constant}{q}
\newcommand{\constanti}{\constant_i}

\newcommand{\ampl}{5}
\newcommand{\freq}{4}

\newcommand{\layers}{2}
\newcommand{\neurons}{300}

\newcommand{\stepvalue}{10^{-4}}

\newcommand{\expa}{a\du{i}}
\newcommand{\expb}{b\du{i}}
\newcommand{\expc}{c\du{i}}

\newcommand{\expterm}{\expa e^{-\expb\T  \footnotesize \begin{bmatrix} 
	x\du{i} \\
	\sigma(x)
\end{bmatrix} + \expc}}

\newcommand{\kk}{\kappa_1}
\newcommand{\kkk}{\kappa_2}

\newcommand{\scalefig}{0.35}

\newcommand{\scaleplot}{0.62}

\def\er/{Erd\H{o}s-R\'enyi}

\def\algoExt/{DEep-Learning aggregative TrAcking}
\def\algo/{DELTA}

\title{Data-Driven Distributed Optimization via Aggregative Tracking and Deep-Learning} 

\author{ 
Riccardo Brumali, \IEEEmembership{Student Member, IEEE},
Guido Carnevale, \IEEEmembership{Member, IEEE},
Giuseppe Notarstefano, \IEEEmembership{Member, IEEE}
	\thanks{Funded by the European Union - NextGenerationEU under the National Recovery and Resilience Plan (PNRR) - Mission 4 Education and research - Component 2 From research to business - Investment 1.1 Notice Prin 2022 -  DD N. 104 del 2/2/2022, from title ECODREAM Energy COmmunity management: DistRibutEd AlgorithMs and toolboxes for efficient and sustainable operations, proposal code 202228CTKY002 - CUP J53D23000560006.
	The authors are with the Department of Electrical,  Electronic and Information Engineering,  Alma Mater Studiorum - Universita` di Bologna,  Bologna, Italy, e-mail: \{name.surname\}@unibo.it.%
}
}

\begin{document}
\maketitle

\begin{abstract}
	In this paper, we propose a novel distributed data-driven optimization scheme.
	In particular, we focus on the so-called aggregative framework, namely, the scenario in which a set of agents aim to cooperatively minimize the sum of local costs, each depending on both local decision variables and an aggregation of all of them.
	We consider a data-driven setup in which each objective function is unknown and can be only sampled at a single point per iteration (thanks to, e.g., feedback from human users or physical sensors).
	We address this scenario through a distributed algorithm that combines three key components: (i) a learning part that leverages neural networks to learn the local cost functions descent direction, (ii) an optimization routine that steers the estimates according to the learned direction to minimize the global cost, and (iii) a tracking mechanism that locally reconstructs the unavailable global quantities.
	By using tools from system theory, i.e., timescale separation and averaging theory, we formally prove that, in strongly convex setups, the overall distributed strategy linearly converges in a neighborhood of the optimal solution whose radius depends on the given accuracy capabilities of the neural networks.
	Finally, we corroborate the theoretical results with numerical simulations.
\end{abstract}


\section{Introduction}

Distributed optimization has recently gained significant attention across various domains.
The comprehensive overview provided in the surveys%
~\cite{giselsson2018large, nedic2018distributed, yang2019survey} discusses the most popular setups and the algorithms used to address these challenges.

The mentioned surveys do not include the so-called aggregative optimization framework which is particularly suited to model tasks arising in cooperative robotics, see the recent tutorial~\cite{testa2023tutorial}.
This novel distributed optimization scenario has been introduced by the pioneering work~\cite{li2021distributed} and deals with networks of agents aiming at cooperatively minimizing the sum of local functions each depending on both microscopic variables (e.g., the position of a single robot) and macroscopic ones (e.g., the spatial distribution of the entire team).
Online and constrained versions of the aggregative problem are considered in%
~\cite{li2021distributed2}.
Other works on distributed aggregative optimization include~\cite{wang2023distributed} in which a distributed method based on the Frank-Wolfe update is proposed, and~\cite{chen2024achieving,wang2024momentum} in which momentum-based algorithms are proposed.
Authors in~\cite{chen2023distributed} propose a continuous time algorithm with non-uniform gradient gains which only requires the sign of relative state information between agents’ neighbors. 
While in~\cite{chen2024compressed} a compressed communication scheme is interlaced with the solution proposed in~\cite{li2021distributed}.
Aggregative optimization scenarios with uncertain environments have been explored in~\cite{carnevale2022learning} where the distributed scheme proposed in~\cite{li2021distributed} is enhanced with a Recursive Least Square (RLS) method estimating the unknown cost via feedback from the users.
A similar framework is addressed in~\cite{brumali2024deep} where instead the learning part relies on a neural network.

Unknown environments as those mentioned above are particularly interesting in the context of the so-called \emph{personalized} optimization frameworks.
In this field, the goal is to minimize cost functions given by the sum of a known part, named engineering function and related to measurable quantities (e.g., time or energy), and an unknown part representing the user's (dis)satisfaction with the current solution.
Since synthetic models based on human preferences often perform well only on average, failing to address the unique preferences of individual users, personalized optimization emphasizes data-driven strategies.
These approaches leverage feedback from specific users about the current solution to adapt and better meet their individual needs.  
With growing human–robot collaboration, personalized optimization techniques can greatly enhance the control of cooperative robotic networks.
In centralized optimization, an initial step toward incorporating user feedback to define human discomfort was undertaken by~\cite{luo2020socially}, where a trajectory design problem was addressed using a cost function influenced by human complaints.  
In~\cite{simonetto2021personalized}, personalized optimization is explored by integrating a learning mechanism based on Gaussian Processes (GP) with an optimization approach. 
Similarly,~\cite{fabiani2022learning} applies a personalized framework within the context of game theory.  
As for personalized distributed frameworks, the articles~\cite{ospina2022time} and~\cite{notarnicola2022distributed} respectively use GP combined with a primal-dual method and RLS merged with the gradient tracking algorithm.

Zeroth-order optimization methods are a popular solution to address problems with unknown cost functions, see~\cite{liu2020primer} for an overview on this topic.
Within this class, the so-called $1$-point optimization methods are particularly appealing since they only require a single function evaluation per iteration and, thus, are particularly suited in scenarios where cost functions evaluations are expensive or multi-point evaluations are not possible, such as when the environment is non-stationary (see~\cite{zhang2022new}).
Examples arise in online optimization, control, reinforcement learning settings (e.g.~\cite{fazel2018global}).
These methods have been widely explored in the centralized literature (see the early references~\cite{flaxman2004online,saha2011improved} or the most recent ones~\cite{chen2022improve,zhang2022new,xiao2023lazy,zhang2024boosting}), while a recent extension to the distributed setting is proposed in~\cite{mhanna2024zero}.
Despite their efficiency, $1$-point methods typically underperform compared to multi-evaluation zeroth-order methods.
As we will better detail later, we overcome this issue by embedding the use of neural networks in our $1$-point architecture.

The use of neural networks in optimization and control is a recent interesting research trend.
For example, in~\cite{cothren2022online,cothren2023perception}, they are integrated into feedback controllers to estimate the system state from perceptual information. 
In~\cite{pirrone2023data}, neural networks are combined with a zeroth-order scheme, while~\cite{schwan2023stability} introduces a general framework for the verification of neural network controllers.
Authors in~\cite{furieri2022neural} employ neural networks to learn stabilizing policies for nonlinear systems.
The work~\cite{martin2024learning} proposes a framework to learn from data algorithms for optimizing nonconvex functions.
Further neural networks have been widely adopted in model predictive control algorithms, where they act as approximators for control policies. 
In this context,~\cite{fabiani2022reliably} proposes a method to certify the reliability of the neural network-based controller, while in~\cite{lupu2024exact} the problem is addressed via HardTanh deep neural networks.
Finally, in~\cite{nubert2020safe} a set point tracking scenario for a robot manipulator is considered.

In this paper, we introduce \algoExt/ (\algo/), a novel distributed data-driven optimization scheme for aggregative problems. 
The peculiarity of \algo/ is to tackle the typical convergence issues of standard $1$-point methods based on static gradient approximations by integrating a dynamic learning strategy based on neural networks that does not increase the number of cost evaluations per iteration.
More in detail, \algo/ operates on a single timescale and integrates three core components: (i) a learning-oriented module, (ii) an optimization-oriented routine, and (iii) a tracking-oriented mechanism.
The learning-oriented component uses local neural networks to asymptotically estimate the correct costs (and their gradients) in a data-driven fashion.
Specifically, each neural network employs only a single cost evaluation per iteration taken in the neighborhood of the current local estimates.
The optimization-oriented component applies an approximated distributed gradient method to the aggregative problem. 
Its inexactness arises from two factors: the fact that the cost functions are unknown and approximated by the networks and the need for global quantities that are locally unavailable.
The tracking-oriented component addresses the second challenge by locally reconstructing the global quantities, i.e., the aggregative variable and the derivative of the global cost with respect to that variable.
We analyze \algo/ using tools from system theory based on timescale separation and averaging theory to formally prove that, in strongly convex settings, the algorithm linearly converges in a neighborhood of the optimal solution whose radius depends on the given accuracy capabilities of the neural networks.
A preliminary version of this work focusing on the time-varying settings is available in~\cite{brumali2024deep}. 
However, in that version, the addressed scenario is significantly simpler. 
Indeed, at each iteration, agents can (i) access multiple samples of their cost functions and (ii) completely train their neural networks until convergence between two consecutive optimization steps.
In contrast, here, at each iteration, we can access the cost functions only at a single point and the neural networks are trained only for a single step.
From a technical point of view, this translates in suitably integrating also the learning dynamics analysis in the optimization and consensus one.
Additionally, in~\cite{brumali2024deep}, the local cost functions are partially known, while in this paper, they are completely unknown.

The paper unfolds as follows. 
In Section~\ref{sec:problem}, we present the problem setup.
In Section~\ref{sec:algo_design}, we introduce the \algo/ algorithm. %
The main theoretical result is provided in Section~\ref{sec:algo_theorem}, while its analysis is detailed in Section~\ref{sec:analysis}.
Finally, in Section~\ref{sec:numerical_simulations}, we  validate our analysis via numerical simulations.

\paragraph*{Notation}

We denote the column stacking of vectors $x_1,\dots,x_N$ with $\col(x_1,\dots,x_N)$.
The $m$-dimensional identity matrix is denoted by $I_m$. %
The symbols $1_N$ and $0_m$ denote the vectors of $N$ ones and $m$ zeros, respectively, while $\1_{N,d} \coloneq 1_N \otimes I_d$, where $\otimes$ denotes the Kronecker product.
Dimensions are omitted when they are clear from the context. 
Given $f: \R^{n_1} \times \R^{n_2} \times \R^{n_3} \to \R^n$, we define $\nabla_1 f(x,y, \theta) \coloneq \dfrac{\partial}{\partial s}f(s,y,\theta)|_{s = x}\T$, $\nabla_2 f(x,y, \theta) \coloneq \dfrac{\partial}{\partial s}f(x,s,\theta)|_{s = y}\T$, and $\nabla_3 f(x,y,\theta) \coloneq \dfrac{\partial}{\partial s}f(x,y,s)|_{s = \theta}\T$.
Given $v \in \R^{n}$ and $r > 0$, we define $\cB_{r}(v) \coloneq \{x \in \R^{n} \mid \norm{x - v} \leq r\}$.
Given the matrices $M_1, \dots, M_N$, we use $\blkdiag(M_1,\dots,M_N)$ to denote the block-diagonal matrix having $M_i$ on the $i$-th block.

\section{Problem Formulation}
\label{sec:problem}

\begin{figure}[hbtp]
	\centering
	\includegraphics[scale=0.66]{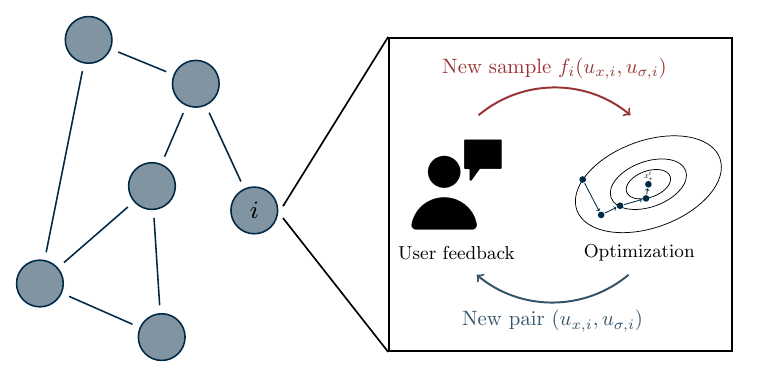}
	\caption{Graphical representation of the problem framework.}
	\label{fig:framework}
\end{figure}

We consider $N$ agents that aim to cooperatively solve
\begin{align}\label{eq:aggregative_optimization_problem}
	\begin{split}
		\min_{(x\du1,\dots,x\du{N}) \in \R^{\n}} \: & \: \sum_{i=1}^{N}\f\du{i}\left(x\du{i},\sigma(x)\right),
	\end{split}
\end{align}
in which $x \coloneq \col(x\du{1}, \dots, x\du{N}) \in \R^\n$ is the global decision vector, with each $x\du{i} \in \R^{\n\du{i}}$ and $\n \coloneq \sum_{i=1}^{N} \n\du{i}$.
For all $i \in \until{N}$, each function $\f\du{i}: \R^{\n\du{i}} \times \R^d \to \R$ represents the \emph{unknown} local objective function of agent $i$ depending on both the local decision variable $x\du{i}$ and the aggregative one given by
\begin{align}\label{eq:sigma}
    \sigma(x) \coloneq \dfrac{1}{N}\sum_{i=1}^{N}\phii(x_i),
\end{align}
where each $\phii: \R^{\n\du{i}} \to \R^d$ is the $i$-th contribution to $\sigma$.
As already mentioned above, the functions $\f\du{i}$ are unknown.
More in detail, we focus on data-driven scenarios where each agent $i$ can only receive a single feedback per iteration about the unknown local cost $\f\du{i}$.
We formalize this aspect as follows.
\begin{assumption}\label{ass:unknown}
	For all $i \in \until{N}$, $\f\du{i}$ is unknown but agent $i$ can receive a single feedback $\f\du{i}(\uxi,\usi)$ per iteration where $(\uxi,\usi) \in \R^{\n\du{i}} \times \R^{d}$ can be arbitrarily chosen.\oprocend
\end{assumption}
The choice of $(\uxi,\usi)$ will represent a design aspect of the overall distributed algorithm we aim to propose.
As one may expect and as we will see later, agent $i$ will choose $(\uxi,\usi)$ in the neighborhood of the current local estimates about (i) the $i$-th block of a problem solution and (ii) the aggregative variable.
The arising data-driven scenario is sketched in Fig.~\ref{fig:framework}.
Now, we formally characterize the functions appearing in problem~\eqref{eq:aggregative_optimization_problem}.
To this end, we first introduce $\fs: \R^{\n} \to \R$ to denote the overall cost function of problem~\eqref{eq:aggregative_optimization_problem}, namely 
\begin{align*}
	\fs(x) \coloneq \sum_{i=1}^N \f\du{i}(x\du{i},\sigma(x)).
\end{align*}
\begin{assumption}\label{ass:objective_functions}
The function $\fs$ is $\mufx$-strongly convex and its gradient is $\Lfs$-Lipschitz continuous, for some $\mufx, \Lfs > 0$.
Further, for all $i \in \until{N}$, $\nabla_1 \f\du{i}$, $\nabla_2 \f\du{i}$, and $\phii$ are Lipschitz continuous with parameters $\Lfx, \Lfsigma, \Lphi > 0$, respectively.
\oprocend
\end{assumption}
We remark that Assumption~\ref{ass:objective_functions} implies the existence of a unique solution $\xstar \in \R^{\n}$ to problem~\eqref{eq:aggregative_optimization_problem}. 

In this paper, we want to develop an optimization algorithm to iteratively solve problem~\eqref{eq:aggregative_optimization_problem} in a distributed manner.
Namely, we want to design an algorithm in which each agent only uses local information and exchanges data with its neighbors.
Indeed, the considered $N$ agents communicate  according to a graph $\cG =(\set,\cE,\adj)$, where
$\set$ is the set of agents, $\cE \subseteq \set \times \set$ is the set of edges, and $\adj \in \R^{N \times N}$ is the weighted adjacency matrix whose
$(i,j)$-entry satisfies $\adjentry_{ij} > 0$ if $(j,i) \in \cE$ and
$\adjentry_{ij} = 0$ otherwise.
The symbol $\cN_i$ denotes the in-neighbor set of agent $i$, namely
$\mathcal{N}_i \coloneq \{j\in\set \mid (j,i)\in\cE\}$.
The next assumption formalizes the class of graphs considered in this work.
\begin{assumption}\label{ass:graph}
	The graph $\cG$ is strongly connected and the matrix $\adj$ is doubly stochastic.\oprocend
\end{assumption}

\section{\algo/: Distributed Algorithm Design}
\label{sec:algo_design}

\begin{figure}%
	\centering
	\includegraphics[scale=0.8]{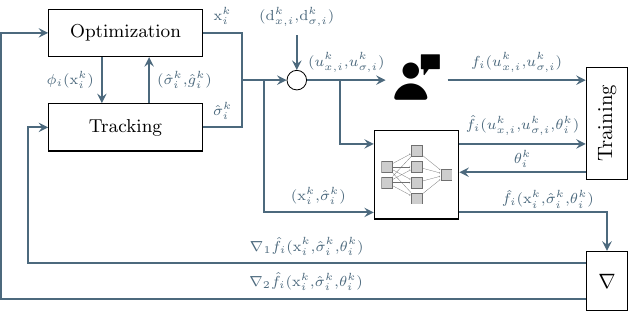}
	\caption{Block diagram of the proposed distributed algorithm, where $\hat{\sigma}_i^\iter \coloneq \wit + \phii(\xit)$ and $\hat{g}_i^\iter \coloneq \zit + \nabla_2 \hfi(\xit,\hat{\sigma}_i^\iter,\thi)$.}
	\label{fig:block}
\end{figure}

In this section, we show the design of \algoExt/ (\algo/), i.e., a novel data-driven distributed method to iteratively address the problem formalized in Section~\ref{sec:problem}.
\algo/ exploits the concurrent action of a learning-oriented part tailored to estimating the unknown gradients of $\f\du{i}$, an optimization-oriented one aimed at solving~\eqref{eq:aggregative_optimization_problem}, and a tracking one devoted to reconstructing the unavailable global quantities.
The combination of these updates results in the \algo/ algorithm, whose graphical description is provided in Fig.~\ref{fig:block} and formal description is reported in Algorithm~\ref{alg:delta}.
\begin{algorithm}%
	\begin{algorithmic}

	\vspace{0.2cm}
	\State \textbf{Learning Update}:
		\begin{align*}
			\thitp = \thit - \step\nabla_3\lossi(\xit + \dxit,\sit + \phii(\xit) + \dsit, \thit)
		\end{align*}
	\State \textbf{Optimization Update}:
	\begin{align*}
		\xitp &= \xit - \step\left[\nabla_1 \hfi\left(\xit,\sit + \phii(\xit), \thit\right) 
		+ \nabla\phii(\xit)\left(\yit + \nabla_2 \hfi\left(\xit,\sit + \phii(\xit), \thit\right)\right)\right] 
	\end{align*}
	\State \textbf{Tracking Update}:
		\begin{align*}
			\sitp &= \sum_{j\in \cN_i}\adjentry_{ij}\left(\sjt + \phij\left(\xjt\right)\right) - \phii\left(\xit\right)
			\\
			\yitp &= \sum_{j\in \cN_i}\adjentry_{ij}\left(\yjt + \nabla_2 \hf\du{j}\left(\xjt, \sjt + \phij\left(\xjt\right), \thjt\right)\right)
			- \nabla_2 \hf\du{i}\left(\xit, \sit + \phii\left(\xit\right), \thit\right)
		\end{align*}

	\vspace{-0.2cm}
	\end{algorithmic}
	\caption{\algo/ (agent $i$)}
	\label{alg:delta}
\end{algorithm}%
Although we analyze \algo/ in a static and deterministic scenario, we point out that its design is particularly well-suited for setups in which the cost functions vary over time and where the cost samples are possibly subjected to disturbances.
Indeed, the synergy between the learning, optimization, and tracking components of \algo/ enables the algorithm to adapt to such changes without requiring a full restart of the neural network training, the optimization task, or the tracking process.
In the following we detail the three main components of \algo/. %

\subsection{Learning Update} 

We introduce in each agent $i$ a local neural network to reconstruct $\f\du{i}$ and the gradients $\nabla_1 \f\du{i}$ and $\nabla_2\f\du{i}$ from data. 
Let $\hfi: \R^{\n\du{i}} \times \R^d \times \R^{\dimthetai} \to \R$ be the estimate of $\f\du{i}$, where $\dimthetai \in \N$ is the number of parameters of the neural network. 
In other words, for a generic parameter $\thi \in \R^{\dimthetai}$, we consider %
\begin{align*}%
	\hfi(x_i,s\du{i},\thi) \proxy \f\du{i}(x\du{i},s\du{i}).
\end{align*}
Once each network structure is fixed (i.e., the number of layers, the number of neurons, and the activation functions are chosen), the function $\hfi$ is analytically available. 
Thus, by considering differentiable functions $\hfi$, we also introduce their gradients%
\begin{subequations}
	\label{eq:approximation}
	\begin{align}
		\nabla_1 \hfi(x_i,s\du{i},\thi)&\proxy \nabla_1\f\du{i}(x_i,s\du{i})
		\\
		\nabla_2 \hfi(x_i,s\du{i},\thi) &\proxy \nabla_2\f\du{i}(x_i,s\du{i}).
	\end{align}
\end{subequations}
We remark that these derivatives are typically available in most of the deep learning packages via automatic differentiation.
The next assumption characterizes the estimates $\nabla_1 \hfi$ and $\nabla_2 \hfi$.
\begin{assumption}\label{ass:hat_lipschitz}
	For all $i \in \until{N}$, the function $\hfi$ is differentiable and its gradients $\nabla_1 \hfi$ and $\nabla_2 \hfi$ are Lipschitz continuous with parameters $\Lhfx, \Lhfsigma > 0$, respectively.\oprocend
\end{assumption}
\begin{remark}
	The requirement in Assumption~\ref{ass:hat_lipschitz} is satisfied by choosing differentiable activation functions with Lipschitz‐continuous gradients, such as tanh, sigmoid, or softplus.
	\oprocend
\end{remark}
\noindent
To iteratively train each neural network using new sampled data, coherently with Assumption~\ref{ass:unknown}, we exploit feedback $\f\du{i}(\xit + \dxit,\sit + \phi(\xit) + \dsit)$, with $\dxit \in \R^{\n\du{i}}$ and $\dsit \in \R^{d}$.
As we will formalize in the next, $\dxit$ and $\dsit$ are two dither signals that we add to force persistency-of-excitation-like properties in the learning process.
Formally, given $\iter \in \N$ and a generic pair $(\xii,\si) \in \R^{\n\du{i}}\times \R^{d}$, we consider the learning problem 
\begin{align}\label{eq:learning_problem}
	\min_{\thi \in \R^{\dimthetai}} \frac{1}{\iter}\sum_{\initer=0}^\iter \lossi(\xii + \dxita,\si+ \dsita, \thi),
\end{align}
where each $\lossi: \R^{\n\du{i}} \times \R^{d} \times \R^{\dimthetai} \to \R$ is the loss function computed using a sample of the observed cost. %
A popular (but not unique) choice for $\lossi$ is the squared loss, namely
\begin{align}
	\lossi(\uxi, \usi, \thi) \!\coloneq\! \frac{1}{2}\norm{\f\du{i}(\uxi, \usi) \!-\! \hfi(\uxi, \usi, \thi)}^2\!.
	 \label{eq:lossi}
\end{align}
The next assumption characterizes the dithers' sequences $\{\dxita,\dsita\}_{\initer\in \N}$ and problem~\eqref{eq:learning_problem}.
\begin{assumption}\label{ass:NN_2}
	There exist $\Tp \in \N \setminus \{0\}$, $\multh > 0$, $\thistar \in \R^{\dimthetai}$, and $\bar{r}\du{i} > 0$
	such that, for all $i \in \until{N}$, the sequence $\{\dxita,\dsita\}_{\initer\in \N}$ is $\Tp$-periodic and, for all $(\xii,\si) \in \R^{\n\du{i}} \times \R^{d}$,  problem~\eqref{eq:learning_problem} with $\iter = \Tp$ is $\multh$-strongly convex in $\thBall$ and $\thistar$ is its local minimum in $\thBall$.
	Further, for all $i \in \until{N}$, $\nabla_3\lossi(\uxi,\usi, \thi)$ is $\Llth$-Lipschitz continuous for all $(\uxi,\usi) \in \R^{\n\du{i}} \times \R^{d}$ and some $\Llth > 0$.\oprocend
\end{assumption}
\begin{remark}\label{rem:ass_2}
	In the ideal case in which $(\xii,\si)$ is fixed and the $\Tp$ samples taken in its neighborhood are all concurrently available, Assumption~\ref{ass:NN_2} ensures (i) local strong convexity of each learning problem~\eqref{eq:learning_problem} into $\thBall$ and (ii) that $\thistar$ is the unique minimizer of its problem~\eqref{eq:learning_problem} into $\thBall$.\oprocend
\end{remark}
Once the local minimizers $\thistar$ have been introduced, we characterize as follows their reconstruction capabilities.
\begin{assumption}\label{ass:NN}
	For all $i \in \until{N}$, there exist Lipschitz-continuous functions $\pertoi : \R^{\n} \to \R^{\n\du{i}}$, $\pertti : \R^{\n} \to \R^{d}$, and
	$\pertzi: \R^{\n} \times \N \to \R^{\dimthetai}$ such that, for all $\x \in \R^\n$ and $\iter \in \until \Tp$, it holds
	\begin{subequations}
		\label{eq:bound_reconstruction_error}
		\begin{align}
			&\nabla_1 \hfi(\xii,\sigma(\x), \thistar)  =  \nabla_1\f\du{i}(\xii,\sigma(\x))  + \pertoi(\x)
			\label{eq:theta_star_nabla_1}
			\\
			&\nabla_2 \hfi(\xii,\sigma(\x),\thistar) =  \nabla_2\f\du{i}(\xii,\sigma(\x)) + \pertti (\x)
			\label{eq:theta_star_nabla_2}
			\\
			&\nabla_3\lossi(\xii + \dxit,\sigma(\x) +\dsit,\thistar) = \pertzi(\x, \iter).
			\label{eq:theta_star_nabla_ell}
		\end{align}
	\end{subequations}
	Moreover, there exists $\epsilon > 0$ such that 
	\begin{align}\label{eq:approx_bound}
		\norm{\pertoi(\x)} &\leq \epsilon, 
		\quad 
		\norm{\pertti(\x)} \leq \epsilon,
		\quad 
		\norm{\pertzi(\x, \iter)} \leq \epsilon,
	\end{align}
	for all $\x \in \R^\n$, $\iter \in \until \Tp$, and $i \in \until{N}$.\oprocend
\end{assumption}
\begin{remark}
	Assumption~\ref{ass:NN} specifies each neural network’s ability to approximate the functions $\f\du{i}$ and their gradients. This capability is captured by the parameter $\epsilon$, the maximum reconstruction error at convergence of the training process (i.e., when $\thi = \thistar$).
	Its value depends on the network design and the unknown $\f\du{i}$, but it is not required by the algorithm.
	\oprocend
\end{remark}%
Once the capabilities of the neural network have been enforced, we can formalize the learning mechanism embedded in \algo/.
An immediate but computationally expensive approach would be to apply the gradient descent method to problem~\eqref{eq:learning_problem}, iteratively updating an estimate $\thit \in \R^{\dimthetai}$ as
\begin{align}\label{eq:plain_GD}
	\thitp \! = \! \thit \!-\! \step \sum_{\initer=0}^{\iter} \nabla_3 \lossi\big(\xit \! + \! \dxita,\wit \! + \! \phii(\xit) \!+ \! \dsita,\thit\big).\!
\end{align}
However, the strategy in~\eqref{eq:plain_GD} would require storing the entire cost samples history and, more importantly, computing an increasing number $\iter$ of gradients in the current estimate $\thit$.
To replace the computationally prohibitive method~\eqref{eq:plain_GD}, we adopt a data-driven update based on a single sample per $\iter$, namely
\begin{align}\label{eq:local_learning_process}
	\thitp = \thit - \step\nabla_3\lossi(\xit + \dxit, \sit + \phii(\xit) + \dsit, \thit).
\end{align}
Indeed, with an eye to the structure formalized in~\eqref{eq:lossi}, we underline that~\eqref{eq:local_learning_process} can be implemented by only relying on a single feedback $\f\du{i}(\xit+\dxit, \sit + \phii(\xit) + \dsit)$ in a neighborhood of $(\xit, \sit + \phii(\xit))$, i.e., of the local estimate of agent $i$ about the $i$-th solution block $\xstari$ and $\sigma(\xt)$, respectively.
Finally, we compactly represent the global cost addressed by the neural networks via $\loss: \R^{\n} \times \R^{Nd} \times \R^{\dimtheta} \to \R$ defined as %
\begin{align}
	\loss(x,s,\theta) \coloneq \sum_{i=1}^N \lossi(x_i,s_i,\theta\du{i}),\label{eq:loss}
\end{align}
where $s \!\coloneq\! \col(s\du{1},\dots,s\du{N}), \theta \!\coloneq\! \col(\theta\du{1},\dots,\theta\du{N})$.
With the notation of Assumption~\ref{ass:NN} at hand, we also define $\dimtheta \coloneq \sum_{i=1}^N \dimthetai$, $\thSet \coloneq \prod_{i=1}^N\thBall \subseteq \R^{\dimtheta}$, and $\thstar \coloneq \col(\thstar\du{1},\dots,\thstar\du{N}) \in \R^{\dimtheta}$.

\subsection{Optimization Update} 

In a full information setup, problem~\eqref{eq:aggregative_optimization_problem} could be solved with a parallel implementation of the gradient method.
Namely, by denoting with $\xit \in \R^{\n\du{i}}$ the estimate of agent $i$ at iteration $\iter$ about the $i$-th block $\xstari \in \R^{\n\du{i}}$ of $\xstar$, each agent would run 
\begin{align}\label{eq:desired_update}
	\xitp &= \xit - \step\nabla_1 \f\du{i}\left(\xit,\sigma(\xt)\right) 
	-\step\frac{\nabla\phii(\xit)}{N}\sum_{j=1}^N \nabla_2 \f\du{j}\left(\xjt,\sigma(\xt)\right),
\end{align}  
where $\step > 0$ is the step size.
Nevertheless, update~\eqref{eq:desired_update} requires the gradients of the unknown $\f\du{i}$. 
We therefore modify it using the approximations~\eqref{eq:approximation}, namely
\begin{align}\label{eq:update_approx_grad}
	\xitp &= \xit - \step\nabla_1 \hf\du{i}\left(\xit,\sigma(\xt),\thit\right) 
	-\step\frac{\nabla\phii(\xit)}{N}\sum_{j=1}^N \nabla_2 \hf\du{j}\left(\xjt,\sigma(\xt), \thjt\right).
\end{align}
Additionally, the unavailable global quantities $\sigma(\xt)$ and $\sum_{j=1}^N \nabla_2 \hf\du{j}(\xjt,\sigma(\xt), \thjt)/N$ are needed in~\eqref{eq:update_approx_grad}.
We thus interlace~\eqref{eq:update_approx_grad} with a tracking-oriented mechanism.

\subsection{Tracking Update}

Being $\sigma(\xt)$ and $\sum_{j=1}^N \nabla_2 \hf\du{j}(\xjt,\sigma(\xt),\thjt)/N$ unavailable, we introduce two auxiliary variables $\sit, \yit \in \R^{d}$ for all $i \in \until{N}$ to estimate them via $\sit + \phii(\xit)$ and $ \nabla_2 \hf\du{i}(\xit, \sit + \phii(\xit),\thit) + \yit$, respectively.
To this end, we update them according to the perturbed consensus dynamics
\begin{subequations}\label{eq:desired_consensus_dynamics}
	\begin{align}
		\sitp &= \sum_{j\in \cN_i}\adjentry_{ij}\left(\sjt + \phij(\xjt)\right) - \phii(\xit)
		\\
		\yitp &= \sum_{j\in \cN_i}\adjentry_{ij}\left(\yjt + \nabla_2 \hf\du{j}(\xjt, \sjt + \phij(\xjt),\thjt)\right)
		- \nabla_2 \hf\du{i}(\xit, \sit + \phii(\xit),\thit), \label{eq:desired_consensus_dynamics_y}
	\end{align}
\end{subequations}
where we recall that the weights $\adjentry_{ij}$ are the entries of the weighted adjacency matrix $\adj$ matching the graph $\cG$.
We remark that the update~\eqref{eq:desired_consensus_dynamics} is fully distributed since it requires only local information and neighboring communication.

\subsection{\algo/: Convergence Properties}
\label{sec:algo_theorem}

To shorten the notation, consider $\tht \!\coloneq\! \col(\tht_1,\dots,\tht_N)$, $\xt \!\coloneq\! \col(\xt_1,\dots,\xt_N)$, $\wt \!\coloneq\! \col(\wt_1,\dots,\wt_N)$, $\qt \!\coloneq\! \col(\qt_1,\dots,\qt_N)$, $\phi(\xt) \!\coloneq\! \col(\phi(\xt_1),\dots,\phi(\xt_N))$, and for all $c \in \{1,2\}$ we define %
\begin{align}
	\hu\du{c}\left(x, w, \theta\right) &\coloneq 
	\begin{bmatrix}
		\nabla\du{c} \hf\du{1}\left(x\du{1},w\du{1}, \theta\du{1}\right) \\
		\vdots \\
		\nabla\du{c} \hf\du{N}\left(x\du{N},w\du{N}, \theta\du{N}\right)
	\end{bmatrix}. 
	\label{eq:ust}
\end{align}
Let the function $\Xi : \R^{\n} \!\times \R^{Nd} \!\times \R^{Nd} \!\times \R^{\dimtheta} \to \R^{\n+2Nd+\dimtheta}$ be
\begin{align*}
	\!\!\Xi(\x, \w, \z, \thh) 
	\!&\coloneq\!\!
	\begin{bmatrix}
		\x - \xstar
		\\
		\w + \phi(\x) - \1\sigma(\x) 
		\\
		\z + (I - \frac{\1\1\T}{N}) \hust(\x, \w + \phi(\x), \thh)
		\\
		\thh - \thstar
	\end{bmatrix}\!.
\end{align*}
Finally, for all $i\in\until{N}$ consider the set 
\begin{align*}
	\cS\du{i} \coloneq \thBall \times \R^{\n\du{i}} \times \{0_d\} \times \{0_d\}.
\end{align*}

\noindent
The next theorem states the convergence properties of \algo/.
\begin{theorem}\label{th:convergence}
	Consider \algo/ and let Assumptions~\ref{ass:unknown}-\ref{ass:NN} hold.
	Further, assume $(\thh\du{i}\ud{0},\x\du{i}\ud{0},\stracker\du{i}\ud{0},\ytracker\du{i}\ud{0}) \in \cS\du{i}$ for all $i \in \until{N}$.
	Then, there exist $\bar{\step}, \bound, \kk > 0$ and $\kkk \in (0,1)$ such that, for all $\step \in (0,\bar{\step})$, it holds 
	\begin{align*} 
		\norm{\xt - \xstar} \leq \kk\kkk^\iter\norm{\Xi(\x\ud0,\w\ud0,\z\ud0,\thh\ud0)} + \epsilon\bound,
	\end{align*} 
	for all $\iter \in \N$.
	\oprocend
\end{theorem}%
The proof of Theorem~\ref{th:convergence} is provided in Section~\ref{sec:analysis}.
It is based on system theory tools aimed at finding a uniform ultimate bound for the trajectories of the dynamical system describing \algo/.
To this end, the proof relies on discrete-time averaging theory (see Appendix~\ref{app:averaging}) and timescale separation.

Besides the term $\bound$ that depends on many problem parameters (e.g., the network connectivity, the Lipschitz continuity constants), we note that the bound in Theorem~\ref{th:convergence} directly depends on the accuracy capabilities of the neural networks (cf. Assumption~\ref{ass:NN_2}) through the parameter $\epsilon$.
In particular, in the case of perfect estimation (i.e., if $\epsilon = 0$, see Assumption~\ref{ass:NN_2}) Theorem~\ref{th:convergence} ensures exact and linear convergence of \algo/ toward the optimal solution $\xstar$.

\section{\algo/: Stability Analysis}
\label{sec:analysis}

In this section, we carry out the stability analysis of \algo/ to prove Theorem~\ref{th:convergence}. %
Our analysis relies on five main steps.
\begin{enumerate}
	\item In Section~\ref{sec:change_of_coordinates}, we reformulate \algo/ and consider a ``nominal'' version in which the ``finite accuracy'' of the neural networks (see Assumption~\ref{ass:NN}) is compensated. 
	We interpret this nominal scheme as a time-varying two-time-scales system, i.e., as the interconnection of two subsystems identified respectively as \emph{slow} and \emph{fast} dynamics.
	This allows us to study the stability properties of the interconnected system by analyzing two auxiliary dynamics called boundary-layer and reduced systems related respectively to the fast and slow dynamics;
	\item In Section~\ref{sec:boundary_analysis}, we show that the origin is an exponentially stable equilibrium of the boundary layer system;  
	\item In Section~\ref{sec:reduced_analysis}, we carry out the reduced system analysis
	in detail, we show exponential stability of the origin via \emph{averaging theory} (see Theorem~\ref{th:bai} in Appendix~\ref{app:averaging});
	\item In Section~\ref{sec:final}, we exploit the results of steps 2) and 3) to conclude the analysis of the nominal system;
	\item Finally in Section~\ref{sec:final}, we return to the original system interpreting it as a perturbed version of the nominal one. 
	This interpretation allows us to use the stability properties of the nominal system (step 4)) to establish those of the original one, thereby providing the proof of Theorem~\ref{th:convergence}.
\end{enumerate}%
Assumptions~\ref{ass:unknown}-\ref{ass:NN} hold true throughout the entire section.

\begin{remark}
	Although the stability analysis of \algo/ is carried out in a static and deterministic environment, it paves the way for extensions to dynamic and uncertain environments subject to, e.g., variations of $\f\du{i}$ over time and noisy cost samples.\oprocend
\end{remark}

\vspace{-0.15cm}
\subsection{\algo/ as a time-varying two-time-scale system}
\label{sec:change_of_coordinates}

By collecting all the local updates in Algorithm~\ref{alg:delta}, we can compactly rewrite \algo/ as
\begin{subequations}
	\label{eq:casual_sys}
	\begin{align}
		\thtp &= \tht - \step \nabla_3 \loss(\xt + \dxt,\wt + \phi(\xt) + \dst, \tht)
		\\		
		\xtp &= \xt - \step \left[\huxt(\xt, \wt + \phi(\xt), \tht) 
		+ \nabla \phi\left(\xt\right) \left(\hust(\xt, \wt + \phi(\xt), \tht) + \qt\right)\right]
		\\
		\wtp &= \adjkron \wt + \left(\adjkron - I_{Nd}\right) \phi(\xt) 
		\\
		\qtp &= \adjkron \qt + \left(\adjkron - I_{Nd}\right) \hust(\xt, \wt + \phi(\xt), \tht),
	\end{align}
\end{subequations}
where $\dxt\coloneq \col(\dr\ud{\iter}\du{x,1},\dots,\dr\ud{\iter}\du{x,N})$, $\dst\coloneq \col(\dr\ud{\iter}\du{\sigma,1},\dots,\dr\ud{\iter}\du{\sigma,N})$, $\adjkron \coloneq \adj \otimes I_d$.
It is worth noting that the subspace
\begin{align*}
	\cS 
	&\coloneq\! \{\left(\theta, \x,\stracker,\ytracker\right) \!\in\! \R\ud{\dimtheta} \!\times\! \R\ud{\n} \!\times\! \R\ud{Nd} \!\times\! \R\ud{Nd} \!\mid\! \1^\top \! \stracker \!=\! 0, \1^\top \! \ytracker \!=\! 0 \},
\end{align*}
is invariant for system~\eqref{eq:casual_sys}.
Thus, we introduce a change of coordinate to isolate the invariant part of the state and define the error coordinates relative to $\thstar$ and $\xstar$. Formally, it is
\begin{align}
	\begin{bmatrix}
		\theta \\\x \\ \stracker \\ \ytracker
	\end{bmatrix}
	\mapsto
	\begin{bmatrix}
		\tth \\ \tx \\ \wm \\ \wort \\ \qm \\ \qort
	\end{bmatrix} 
	\coloneq
	\begin{bmatrix}
		\theta - \thstar
		\\ 
		\x - \xstar 
		\\ 
		\frac{\1\T}{N}\stracker
		\\
		\Rmat\T\stracker 
		\\
		\frac{\1\T}{N}\ytracker
		\\
		\Rmat\T\ytracker
	\end{bmatrix},\label{eq:change}
\end{align}
with $\Rmat \in \R\ud{Nd \times (N-1)d}$ such that $\Rmat\T\1=0$ and $\Rmat\T\Rmat=I_{Nd}$.  
By using the definitions of $\wmt$ and $\qmt$ in~\eqref{eq:change} and since $\adj$ is doubly stochastic (cf. Assumption~\ref{ass:graph}), it holds 
\begin{align*}
	\wmtp &= \wmt, \quad \qmtp = \qmt,
\end{align*}
which implies $(\wmt,\qmt) = (\winitm,\qinitm)$ for all $\iter \in \N$.
Further, the initialization $\wit = \qit =0$ leads to $(\winitm,\qinitm) = (0,0)$ and, thus, it holds $(\wmt,\qmt) = (0,0)$ for all $\iter \in \N$.
Hence, we rewrite~\eqref{eq:casual_sys} in the new coordinates~\eqref{eq:change} and neglect the trivial dynamics of $\wmt$ and $\qmt$, thus obtaining an equivalent restricted system
\begin{subequations}\label{eq:restricted_system_without_perp}
	\begin{align}
		\tthtp &\!= \ttht - \step \nablalosst\left(\txt,\wortt, \ttht,\iter\right)
		\label{eq:restricted_system_without_perp_tth}
		\\
		\txtp &\!= \txt - \step  \left[\huxt(\txt + \xstar,\Rmat \wortt + \phi(\txt + \xstar), \ttht + \thstar) \right.\notag
		\\
		 &\hspace{9pt} 
		 + \! \nabla \phi(\txt \! + \! \xstar) \hust(\txt \! + \! \xstar\!,\! \Rmat \wortt \! + \! \phi(\txt \! + \! \xstar), \ttht \! + \! \thstar)
		 \notag\\
		 &\hspace{9pt} 
		 + \nabla \phi(\txt + \xstar) \Rmat \qortt \Big]
		\label{eq:restricted_system_without_perp_tx}
		\\
		\worttp &\!= \!\Rmat\T\!\adjkron\Rmat \wortt \!+\! \adjminusI \phi(\txt \!+\! \xstar)\label{eq:restricted_system_without_perp_w}
		\\
		\qorttp &\!=\! \Rmat\T\!\adjkron\Rmat \qortt  \!+\!  \adjminusI \hust(\txt \!+\! \xstar,\wortt \!+\! \phi(\txt \!+\! \xstar), \ttht \!+\! \thstar).\!\!
		\label{eq:restricted_system_without_perp_q}
	\end{align}
\end{subequations}
where $\adjminusI \coloneq \Rmat^\top (\adjkron - I_{Nd})$ and we introduce the compact notation $\tilde{\loss}: \mathbb{R}^{n} \times \mathbb{R}^{(N-1)d} \times \mathbb{R}^{\dimtheta} \times \mathbb{N} \to \mathbb{R}$, defined as
\begin{align*}
	\tilde{\loss}(\tx,w_\perp, \tth, \iter) \coloneq \loss\left(\tx + \xstar + \dxt,\Rmat w_\perp + \phi(x) + \dst, \tth + \thstar\right).
\end{align*}
Now, we introduce the ``nominal'' version of~\eqref{eq:restricted_system_without_perp}, i.e., an auxiliary system in which the finite accuracy of the neural networks (see Assumption~\ref{ass:NN}) is suitably compensated by adding some fictitious terms. %
Such a nominal system reads as
\begin{subequations}\label{eq:nominal_restricted_system}
	\begin{align}
		\tthtp &\!= \ttht \!-\! \step \!\left[\nablalosst\!\left(\txt,\wortt, \ttht,\iter\right)- \pertz\!\left(\txt + \xstar, \iter\right)\right]\!\!
		\label{eq:nominal_restricted_system_tth}
		\\
		\txtp \!&=\! \txt \!-\! \step \! \left[\huxt\!\left(\txt \!+\! \xstar\!,\Rmat \wortt \!+\! \phi\!\left(\txt \!+\! \xstar\right)\!, \ttht \!+\! \thstar\right) \right.\notag\\
		 &\hspace{9pt} 
		 + \!\nabla \phi\!\left(\txt \!\!+\! \xstar\right) \!\hust\!\left(\!\txt \!\!+\! \xstar\!,\Rmat \wortt \!+\! \phi\!\left(\txt \!\!+\! \xstar\right)\!, \ttht \!+\! \thstar\!\right) 
		 \notag\\
		 &\hspace{9pt} %
			+ \!\nabla \phi\!\left(\txt \!+\! \xstar\right)\!\Rmat \qortt 
		\notag\\
		&\hspace{9pt} %
			\!-\! \perto\!\left(\txt \!+\! \xstar\right) \!-\! \nabla \phi\!\left(\txt \!+\! \xstar\right)\! \tfrac{\1\1\T}{N} \pertt\!\left(\txt \!+\! \xstar\right)\! \Big]\!\!
		 \label{eq:nominal_restricted_system_tx}
		\\
		\worttp &\!=\! \Rmat\T\!\adjkron\Rmat \wortt \!+\! \adjminusI \phi\left(\txt \!+\! \xstar\right)\label{eq:nominal_restricted_system_w}
		\\
		\qorttp &\!=\! \Rmat\T\!\adjkron\Rmat \qortt \! + \! \adjminusI \hust\!\left(\txt \!+\! \xstar\!,\wortt \!+\! \phi(\txt \!+\! \xstar), \ttht \!+\! \thstar\right)\!\!,\!\!
		\label{eq:nominal_restricted_system_q}
	\end{align}
\end{subequations}
where $\pert\du{c}(x)\!\coloneq\! \col(\pert\du{c,1}(x), \dots, \pert\du{c,N}(x))$ for all $c \!\in\! \{1,2\}$, while $\pertz\left(x,\iter\right) \!\coloneq\! \col(\pert\du{\loss,1}\!\left(x,\iter\right),\dots,\pert\du{\loss,N}\!\left(x,\iter\right))$.
Due to $\perto$, $\pertt$, and $\pertz$, $(0,0,-\Rmat^\top \phi(\xstar), -\Rmat^\top \hust(\xstar, \mathbf{1}\sigma(\xstar), \thstar))$ is an equilibrium of~\eqref{eq:nominal_restricted_system} (cf. Assumption~\ref{ass:NN}). 
In detail, the equilibrium points of the subsystems~\eqref{eq:nominal_restricted_system_w}-\eqref{eq:nominal_restricted_system_q} are parametrized in the states $(\theta,\x)$ of~\eqref{eq:nominal_restricted_system_tth}-\eqref{eq:nominal_restricted_system_tx} through $\h: \mathbb{R}^{\dimtheta} \times \mathbb{R}^n \to \mathbb{R}^{2(N-1)d}$ defined as
\begin{align}
	\label{eq:parametrized_equilibria}
	\h\left(\theta, x\right)
	 \coloneq 
	\begin{bmatrix}
		- \Rmat\T\phi\left(x\right)\\
		- \Rmat\T\hust\left(x,\1\sigma\left(x\right), \theta\right)
	\end{bmatrix}.
\end{align}
Further, the variations of~\eqref{eq:nominal_restricted_system_tth}-\eqref{eq:nominal_restricted_system_tx} can be arbitrarily reduced via $\step$.
From these observations, we interpret system~\eqref{eq:nominal_restricted_system} as a time-varying two-time-scale system, where~\eqref{eq:nominal_restricted_system_tth}-\eqref{eq:nominal_restricted_system_tx} is the so-called slow dynamics, while~\eqref{eq:nominal_restricted_system_w}-\eqref{eq:nominal_restricted_system_q} is the so-called fast part.
According to this view, we now proceed by analyzing the boundary-layer and reduced system associated to~\eqref{eq:nominal_restricted_system}.

\subsection{Boundary-Layer System Analysis}
\label{sec:boundary_analysis}

In this section, we analyze the boundary-layer system associated with~\eqref{eq:nominal_restricted_system}, obtained by fixing an arbitrary $(\tth, \tx) \in \mathbb{R}^{\dimtheta} \times \mathbb{R}^{n}$ in~\eqref{eq:nominal_restricted_system_w}-\eqref{eq:nominal_restricted_system_q}, written in the error coordinates.
\begin{align*}
	\begin{bmatrix}
		\twortt
		\\ 
		\tqortt
	\end{bmatrix}
	\coloneq
	\begin{bmatrix}
		\wortt
		\\
		\qortt
	\end{bmatrix}
	-
	\h\left(\tth + \thstar,\tx + \xstar\right).
\end{align*}
In light of the stochasticity property of $\adj$ (cf. Assumption~\ref{ass:graph}) and the definition of $\h$ (cf.~\eqref{eq:parametrized_equilibria}), we use the properties $\Rmat\T\1 = 0$ and $\Rmat\Rmat\T = I_{Nd} - \frac{1}{N}\1\1\T$ to write this system as
\begin{subequations}
	\label{eq:boundary_sys}
	\begin{align}
		\tworttp &= \Rmat\T\adjkron\Rmat \twortt 
		\\
		\tqorttp &= \Rmat\T\adjkron\Rmat \tqortt \!+\! \adjminusI \Delta\hust\!\left(\tx \!+\! \xstar, \twortt, \tth \!+\! \thstar\right),
	\end{align}
\end{subequations}
where $\Delta\hust: \R^{n} \times \R^{(N-1)d} \times \R^{\dimtheta} \to \R$ is defined as
\begin{align}
	\Delta\hust(\tx + \xstar, \twort, \tth + \thstar) &\coloneq \hust(\tx + \xstar, \Rmat\twort + \1\sigma(\tx + \xstar), \tth + \thstar) - \hust(\tx + \xstar,\1\sigma(\tx + \xstar), \tth + \thstar).\label{eq::u2}
\end{align}
Let us introduce $\stateF \coloneq \col \left(\stateF_1,\stateF_2\right) \coloneq \col \left(\twort, \tqort\right)$, $\RARdiag \coloneq \blkdiag(\Rmat\T \adjkron \Rmat,\Rmat\T \adjkron \Rmat)$, and $\gFast: \R^{n} \times \R^{2(N-1)d} \times \R^{\dimtheta} \to \R^{2(N-1)d}$ defined as
	\begin{align}
		\gFast\left(\tx, \stateF, \tth\right) &\coloneq 
		\begin{bmatrix}
			0\\
			\Rmat\T\adjkron \Delta \hust (\tx + \xstar, \stateF_1, \tth + \thstar)
		\end{bmatrix},
		\label{eq:gFast}
	\end{align}
		which allow us to compactly rewrite~\eqref{eq:boundary_sys} as 
		\begin{align}
			\stateFtp = \RARdiag\stateFt + \gFast(\tx,\stateFt, \tth).\label{eq:boundary_sys_compact}
		\end{align}
With this notation at hand, we are ready to formalize the stability properties of the boundary-layer system~\eqref{eq:boundary_sys_compact}.
\begin{lemma}\label{lem:boundary_lyap}
	There exist a Lyapunov function $U:\R^{2(N-1)d} \to \R$ and $\bb_1,\bb_2,\bb_3,\bb_4>0$ such that
	\begin{subequations}
		\begin{align}
			&\bb_1 \!\norm{\stateF}^2 \!\leq \Vfast\!\left(\stateF\right) \!\leq \bb_2\! \norm{\stateF}^2 \label{eq:positive_Vfast}\\ 
			&\Vfast\!\left(\!\RARdiag \stateF \!+\! \gFast\!\left(\tx, \stateF, \tth\right)\right) \!-\! \Vfast\!\left(\stateF\right)\leq - \bb_3 \!\norm{\stateF}^2 \label{eq:delta_Vfast_negative}\\
			&\left\lvert \Vfast\!\left(\stateF\right) \!-\! \Vfast\!\left(\bar{\stateF}\right) \right\rvert \leq \bb_4 \!\norm{\stateF - \bar{\stateF}} \left(\norm{\stateF} \!+\! \norm{\bar{\stateF}}\right), \label{eq:lipschitz_like_bound_Vfast}
		\end{align}
	\end{subequations}
	for any $\stateF, \bar{\stateF} \in \R^{2(N-1)d}$ and $\left(\tth, \tx\right)\in \R^{\dimtheta} \times \R^n$. 
	\oprocend
\end{lemma}
\noindent The proof of Lemma~\ref{lem:boundary_lyap} is reported in Appendix~\ref{sec:proof_bl}.

\subsection{Reduced system analysis}
\label{sec:reduced_analysis}

In this section, we analyze the reduced system associated to~\eqref{eq:nominal_restricted_system}, i.e., the system obtained by plugging $\col(\wortt,\qortt) = \h(\ttht + \thstar, \txt + \xstar)$ for all $\iter \in \N$ into~\eqref{eq:nominal_restricted_system_tth}-\eqref{eq:nominal_restricted_system_tx}, namely%
\begin{subequations}
	\label{eq:reduced_sys_imp}
	\begin{align}
		\tthtp \!&= \ttht \!-\! \step (\nablalosst(\!\txt,\1\sigma\!\left(\txt \!+\! \xstar\right), \ttht,\iter) \!-\! \pertz\!\left(\iter\right)) \label{eq:reduced_sys_imp_tth}
		\\
		\!\txtp \!&= \txt \!-\! \step \hud(\txt,\ttht), \label{eq:reduced_sys_imp_tx}
	\end{align}
\end{subequations}
in which we introduce $\hud: \R^{\n} \times \R^{\dimtheta} \to \R^{\n}$ defined as
\begin{align}
	\hud(\tx,\tth) &\coloneq \huxt(\tx + \xstar, \1\sigma(\tx + \xstar), \tth + \thstar) - \perto(\tx + \xstar)
	\notag\\ 
	&\hspace{13pt}+ \nabla\phi(\tx + \xstar) \tfrac{\1\1\T}{N}\hust(\tx + \xstar,\1\sigma(\tx + \xstar), \tth + \thstar)
	\notag\\ 
	&\hspace{13pt}-\nabla\phi\left(\tx + \xstar\right) \tfrac{\1\1\T}{N}\pertt(\tx+\xstar)
	.\label{eq:uu}
\end{align}
Using the accuracy characterization~\eqref{eq:theta_star_nabla_ell} of $\thstar$ and the optimality of $\xstar$ for~\eqref{eq:aggregative_optimization_problem}, we point out that the origin is an equilibrium of~\eqref{eq:reduced_sys_imp}. 
We then study its stability via the corresponding averaged system, obtained by averaging the time-varying vector field of~\eqref{eq:reduced_sys_imp} over an infinite time horizon (see Appendix~\ref{app:averaging}).
To this end, let $\lossav: \R^{\n} \times \R^{\dimtheta} \to \R$ be 
\begin{align*}
	\lossav(\tilde{x},\tilde{\theta}) \coloneq\tfrac{1}{\Tp} \textstyle\sum_{\iter=1}^{\Tp} \loss(\tilde{x} \!+\! \xstar \!+\! \dxt,\1\sigma(\tilde{x} \!+\! \xstar) \!+\! \dst,\tth \!+\! \thstar).
\end{align*}
The periodicity of $\dxt$ and $\dst$ (cf. Assumption~\ref{ass:NN}) ensures the existence of the averaged system associated to~\eqref{eq:reduced_sys_imp_tth}, namely
\begin{subequations}
	\label{eq:reduced_sys_av_with_p}
	\begin{align}
		\tthtpav &\!=\! \tthtav \!-\! \step\nablalossav(\txav,\tthtav) \! + \! \frac{\step}{\Tp}\!\sum_{\iter=1}^{\Tp} \pertz\!\left(\txav \! + \! \xstar\!,\iter\right)\!
		\\
		\!\txtpav \!&= \txtav \!-\! \step \hud(\txtav,\tthtav).
	\end{align}
\end{subequations}
From Assumption~\ref{ass:NN_2} $\nablalossav(\cdot,\tthav)$ is the gradient of a $\multh$-strongly convex function for all $(\tthav + \thstar) \in \thSet$.
We thus obtain the following stability properties of the averaged system~\eqref{eq:reduced_sys_av_with_p}.
\begin{lemma}
	\label{lemma:av}
	There exists $\bar{\step}_2 > 0$ such that, for all $\step \in (0,\bar{\step}_2)$, the origin is an exponentially stable equilibrium of~\eqref{eq:reduced_sys_av_with_p} with domain of attraction $\{(\tthav,\txav) \in \R^{\dimtheta} \times \R^{\n} \mid (\tthav + \thstar) \in \thSet\}$.
\end{lemma}
\noindent The proof of Lemma~\ref{lemma:av} is reported in Appendix~\ref{sec:proof_av}.\\

Let us introduce $\stateS \coloneq \col(\stateS_1, \stateS_2) \coloneq \col(\tth, \tx)$ and $\gSlow: \R^{\dimtheta + \n} \times \N \to \R^{\dimtheta + \n}$ defined as
\begin{align*}
	\gSlow\left(\stateS, \iter\right) &\coloneq
	\begin{bmatrix}
		\gSlowth\left(\stateS_1, \stateS_2, \iter\right)	
		\\
		\gSlowx\left(\stateS_1, \stateS_2 \right)
	\end{bmatrix}
	\\
	&\coloneq
	\begin{bmatrix}
		\nablalosst\left(\stateS_2, \1 \sigma(\stateS_2 + \xstar),\stateS_1, \iter\right) -\pertz(\state_2 + \xstar,\iter)
		\\
		\hud(\stateS_2, \stateS_1)
	\end{bmatrix}\!.
\end{align*}
Then, we can compactly rewrite system~\eqref{eq:reduced_sys_imp} as 
\begin{align}
	\stateStp = \stateSt - \step \gSlow\left(\stateSt, \iter\right).\label{eq:reduced_sys_compact}
\end{align}
Let $\cX \coloneq \{\stateS \in \R^{\dimtheta + \n}\mid  \stateS_1 + \thstar \in \thSet\}$.
The next lemma ensures exponential stability properties of the origin for~\eqref{eq:reduced_sys_compact}.
\begin{lemma}
	\label{lem:reduced_stability}
	There exist a Lyapunov function $\Vconv:\R^{\dimtheta + \n} \times N \to \R$ and $\bar{\step}_1 > 0$ such that, for all $\step \in (0,\bar{\step}_1)$ it holds
	\begin{subequations}
		\begin{align}
			&\cc_1 \!\norm{\stateS}^2 \!\leq \Vconv\!\left(\stateS, \iter\right) \!\leq \cc_2\! \norm{\stateS}^2 \label{eq:positive_Vslow}\\ 
			&\Vconv\!\left(\stateS - \step \gSlow(\stateS, \iter), \iterp\right) \!-\! \Vconv\!\left(\stateS, \iter\right) \leq - \step\cc_3 \!\norm{\stateS}^2 \label{eq:delta_Vslow_negative}\\
			&\left\lvert \Vconv\!\left(\stateS, \iter\right) \!-\! \Vconv\!\left(\bar{\stateS}, \iter\right) \right\rvert \leq \cc_4 \!\norm{\stateS - \bar{\stateS}} \left(\norm{\stateS} \!+\! \norm{\bar{\stateS}}\right), \label{eq:lipschitz_like_bound_Vslow}
		\end{align}
	\end{subequations}
	for all $\stateS, \bar{\stateS} \in \cX$, and some $\cc_1,\cc_2,\cc_3,\cc_4>0$. \oprocend
\end{lemma}
\noindent The proof of Lemma~\ref{lem:reduced_stability} is reported in Appendix~\ref{sec:proof_reduced_stability}.

\subsection{Stability of Nominal System} 
\label{sec:nominal} 
In this section, we combine the stability properties of the boundary-layer (cf. Lemma~\ref{lem:boundary_lyap}) and reduced system (cf. Lemma~\ref{lem:reduced_stability}) to formalize those of the nominal system~\eqref{eq:nominal_restricted_system}.
To this end, let $\nstate \coloneq n + \dimtheta + 2(N-1)d$ and define $\state \in \R^{\nstate}$ as
\begin{align}
	\state \coloneq 
	\begin{bmatrix}
		\tth \\ \tx \\ \col(
			\wort,\qort)- h(\tth + \thstar,\tx + \xstar) 
	\end{bmatrix}
	.\label{eq:state}
\end{align}
With this notation at hand, we rewrite~\eqref{eq:nominal_restricted_system} as 
\begin{align}
	\statetp = \dyn(\statet,\iter),\label{eq:nominal_system_error}
\end{align}
where $\dyn: \R^{\nstate} \times \N \to \R^{\nstate}$ compactly represents the update~\eqref{eq:nominal_restricted_system} in the new coordinates~\eqref{eq:state}.
The next lemma ensures that the origin is an exponentially stable equilibrium point of~\eqref{eq:nominal_system_error}.
We first introduce the set $\domVpert \subset \R^{\nstate}$ defined as  
\begin{align*}
	\domVpert \coloneq \{\state \coloneq \col(\state_1,\state_2,\state_3,\state_4) \in \R^{\nstate} \!\mid\!  \state_1 + \thstar \in \thSet\}.
\end{align*}
\begin{lemma}\label{lemma:nominal}
	There exist a function $\Vnom:\R^{\nstate} \times N \to \R$ and $\bar{\step}, \aa_1,\aa_2,\aa_3,\aa_4 > 0$ such that, for all $\step \in (0,\bar{\step})$, it holds
	\begin{subequations}
		\begin{align}
			&\aa_1 \!\norm{\state}^2 \!\leq \Vnom\!\left(\state, \iter\right) \!\leq \aa_2\! \norm{\state}^2 \label{eq:positive_Vnom}\\ 
			&\Vnom\!\left(\dyn(\state,\iter), \iterp\right) \!-\! \Vnom\!\left(\state, \iter\right) \leq - \step\aa_3 \!\norm{\state}^2 \label{eq:delta_Vnom_negative}\\
			&\left\lvert \Vnom\!\left(\state, \iter\right) \!-\! \Vnom\!\left(\bar{\state}, \iter\right) \right\rvert \leq \aa_4 \!\norm{\state - \bar{\state}} \left(\norm{\state} \!+\! \norm{\bar{\state}}\right), \label{eq:lipschitz_like_bound_Vnom}
		\end{align}
	\end{subequations}
	for all $\state, \bar{\state} \in \domVpert$ and $\iter \in \N$.\oprocend
\end{lemma}
\noindent The proof of Lemma~\ref{lemma:nominal} is reported in Appendix~\ref{sec:proof_nominal}.

\subsection{Proof of Theorem~\ref{th:convergence}}
\label{sec:final}

The proof consists in finding a uniform ultimate bound to the trajectories of system~\eqref{eq:restricted_system_without_perp} by using the Lyapunov function $\Vnom$ studied in Lemma~\ref{lemma:nominal}. %
To this end, by following the notation in~\eqref{eq:nominal_system_error}, we compactly rewrite the original system~\eqref{eq:restricted_system_without_perp} as
\begin{align} \label{eq:original_system}
	\statetp \coloneq \dyno (\statet, \iter) = \dyn(\statet, \iter) - \step\pert(\statet),
\end{align}
with $\dyno: \R^{\nstate} \to \R^{\nstate}$ and $\pert: \R^{\nstate} \to \R^{\nstate}$ which reads as  
\begin{align*} %
	\pert(\state) \coloneq \begin{bmatrix}
		\pertz(\state_2+ \xstar, \iter)  
		\\
		\perto(\state_2 + \xstar) + \nabla \phi \left(\state_2 + \xstar\right)\tfrac{\1\1\T}{N} \pertt(\state_2 + \xstar)
		\\
		\tfrac{\Rmat\T}{\gamma} \left(\phi\left(\dynocomp{2}(\state, \iter) + \xstar\right) - \phi\left(\dyncomp{2}(\state, \iter) + \xstar\right) \right) 
		\\
		\tfrac{\Rmat\T}{\gamma} \left(\hust^{\tny{or}}\left(\dyno\left(\state, \iter\right)\right) -  \hust^{\tny{nom}}\left(\dyn\left(\state, \iter\right)\right)\right)
	\end{bmatrix},
\end{align*}
in which, $\dyncomp{i}(\statet,\iter)$ and $\dynocomp{i}(\statet,\iter)$ denote the $i$-th component of $\dyn(\statet,\iter)$ and $\dyno(\statet,\iter)$, while $\hust^{\tny{or}}(\dyno(\state, \iter))$ and $\hust^{\tny{nom}}(\dyn(\state, \iter))$ represent the function $\hust$ evaluated respectively at $\left(\dynocomp{2}(\state, \iter) \!+\! \xstar, \1 \sigma\left(\dynocomp{2}(\state, \iter) \!+\! \xstar\right), \dynocomp{1}(\state, \iter) \!+\! \thstar\right)$ and $\left(\dyncomp{2}(\state, \iter) \!+\! \xstar, \1 \sigma\left(\dyncomp{2}(\state, \iter) \!+\! \xstar\right), \dyncomp{1}(\state, \iter) \!+\! \thstar\right)$.

We proceed bounding $\norm{\dynocomp{2}(\statet\!,\iter) \!-\! \dyncomp{2}(\statet\!,\iter)}^2$ as
\begin{align}
	\label{eq:bound_xor_xnom}
	\norm{\dynocomp{2}(\statet\!,\iter) \!-\! \dyncomp{2}(\statet\!,\iter)}^2
	&= \big\lVert \perto\!\left(\state_2 \!+\! \xstar\right) \!+\! \nabla \phi \left(\state_2 \!+\! \xstar\right)\tfrac{\1\1\T}{N} \pertt\!\left(\state_2 \!+\! \xstar\right) \big\rVert^2
	\notag
	\\
	&\stackrel{(a)}{\leq} \norm{\perto\left(\state_2 + \xstar\right)}^2 + \norm{\nabla \phi \left(\state_2 + \xstar\right)\tfrac{\1\1\T}{N} \pertt\left(\state_2 + \xstar\right)}^2
	\notag\\
	&\hspace{11pt} +2 \perto\left(\state_2 + \xstar\right)\T \nabla \phi \left(\state_2 + \xstar\right) \tfrac{\1\1\T}{N} \pertt \left(\state_2 + \xstar\right)
	\notag\\
	&\stackrel{(b)}{\leq} \epsilon^2 N \left( 1 + \Lphi^2 + 2 \Lphi \right) \notag
	\\
	&= \epsilon^2 N \left( 1 + \Lphi\right)^2,
\end{align}
where in $(a)$ we expanded the square norm, while in $(b)$ we used the Cauchy-Schwarz inequality, the Lipschitz continuity of $\phi$ (see Assumption~\ref{ass:objective_functions}), and the accuracy bounds enforced by Assumption~\ref{ass:NN}. 
Similarly, we compute the following bound 
\begin{align}
	\label{eq:bound_sigmaor_sigmanom}
	\norm{\1\sigma(\dynocomp{2}(\statet,\iter) + \xstar) - \1\sigma(\dyncomp{2}(\statet,\iter) + \xstar)}^2
	&=
	\norm{\tfrac{\1\1\T}{N} \left(\phi(\dynocomp{2}(\statet,\iter)) - \phi(\dyncomp{2}(\statet,\iter))\right)}^2
	\notag
	\\
	&\leq \Lphi^2 \norm{\dynocomp{2}(\statet,\iter) - \dyncomp{2}(\statet,\iter)}^2 
	\notag\\
	&\stackrel{(a)}{\leq} \epsilon^2 N \Lphi^2 \left( 1 + \Lphi\right)^2, 
\end{align} 
where in $(a)$ we used the bound~\eqref{eq:bound_xor_xnom}.
Lastly, we notice that
\begin{align}
	\label{eq:bound_thor_thnom}
	\norm{\dynocomp{1}(\statet,\iter) - \dyncomp{1}(\statet,\iter)}^2 = \norm{\gamma \pertz(\state_2 + \xstar, \iter)}^2 
	\leq \gamma^2 \epsilon^2 N. 
\end{align} 
At this point, $\norm{\pert(\state)}$ can be bounded as
\begin{align}
	\norm{\pert(\state)}^2 
	&
	= \norm{\pertz(\state_2 \!+\! \xstar, \iter)}^2
	\notag\\
	&\hspace{13pt} \!+\! \norm{\perto\left(\state_2 \!+\! \xstar\right) \!+\! \nabla \phi \!\left(\state_2 \!+\! \xstar\right)\tfrac{\1\1\T}{N} \pertt\left(\state_2 \!+\! \xstar\right)}^2 
	\notag\\
	&\hspace{13pt} \!+\! \norm{\tfrac{\Rmat\T}{\gamma} \!\left(\phi\!\left(\dynocomp{2}(\state\!, \iter) \!+\! \xstar\right) \!-\! \phi\!\left(\dyncomp{2}(\state\!, \iter) \!+\! \xstar\right)\right) }^2
	\notag\\
	&\hspace{13pt} \!+\! \norm{\tfrac{\Rmat\T}{\gamma}\hust^{\tny{or}}\left(\dyno\left(\state, \iter\right)\right) -  \hust^{\tny{nom}}\left(\dyn\left(\state, \iter\right)\right)}^2
	\notag
	\\
	&\stackrel{(a)}{\leq} \norm{\pertz(\state_2 \!+\! \xstar, \iter)}^2
	\notag\\
	&\hspace{13pt} \!+\! \norm{\perto\left(\state_2 \!+\! \xstar\right) \!+\! \nabla \phi \!\left(\state_2 \!+\! \xstar\right)\tfrac{\1\1\T}{N} \pertt\left(\state_2 \!+\! \xstar\right)}^2 
	\notag\\
	&\hspace{13pt} \!+\! \tfrac{\Lphi^2}{\gamma^2} \norm{\Rmat}^2 \norm{\dynocomp{2}(\statet\!,\iter) \!-\! \dyncomp{2}(\statet\!,\iter)}^2
	\notag\\
	&\hspace{13pt} \!+\! \tfrac{\Lfsigma^2}{\gamma^2} \norm{\Rmat}^2 \left( \norm{\dynocomp{2}(\statet\!,\iter) \!-\! \dyncomp{2}(\statet\!,\iter)}^2 \right. 
	\notag
	\\
	&\hspace{13pt} \!+\! \norm{\1\sigma(\dynocomp{2}(\statet\!,\iter) \!+\! \xstar) \!-\! \1\sigma(\dyncomp{2}(\statet\!,\iter) \!+\! \xstar)}^2 
	\notag
	\\
	&\left. \hspace{13pt} \!+\! \norm{\dynocomp{1}(\statet\!,\iter) \!-\! \dyncomp{1}(\statet\!,\iter)}^2\right)
	\notag
	\\
	&\stackrel{(b)}{\leq} \epsilon^2 N \!\Big[\!1 \!+\! \left( 1 \!+\! \Lphi\right)^2 \!+\! \tfrac{\Lphi^2}{\gamma^2} \norm{\Rmat}^2 \left( 1 \!+\! \Lphi\right)^2 
	\notag\\
	&\hspace{13pt} \!+\! \tfrac{\Lfsigma^2}{\gamma^2} \!\norm{\Rmat}^2 \! \left(
		\! \left( 1 \!+\! \Lphi\right)^2\!
		\!+\!
		 \Lphi^2 \left( 1 \!+\! \Lphi\right)^2\!
		\!+\!
		\gamma^2
	\right)\!\!\Big] 
	\notag\\
	& \eqcolon \epsilon^2 \aa_5^2,
	\notag
\end{align}
where in $(a)$ we used the Cauchy-Schwarz inequality and the Lipschitz continuity of $\phii(\xii)$ and $\nabla_2 \hfi$ (cf. Assumptions~\ref{ass:objective_functions} and~\ref{ass:hat_lipschitz}), while in $(b)$ we used equations~\eqref{eq:bound_xor_xnom}-\eqref{eq:bound_thor_thnom}.
With this bound at hand, we now consider the Lyapunov function $\Vnom$ introduced in Lemma~\ref{lemma:nominal} to study the stability properties of the original system~\eqref{eq:original_system}.
Then, by evaluating the increment of $\Vnom$ along the trajectories of~\eqref{eq:original_system}, for all $\iter \in \N$ and $\state \in \domVpert$, we get
\begin{align}
	&\dVnom\left(\state, \iter\right) 
	\coloneq 
	\Vnom\left(\dyn\left(\state, \iter\right) - \step\pert\left(\state\right), \iterp\right) - \Vnom\left(\state, \iter\right) 
	\notag
	\\
	&\stackrel{(a)}{=} \Vnom\left(\dyn\left(\state, \iter\right), \iterp\right) - \Vnom\left(\state, \iter\right) 
	\notag\\
	&\hspace{13pt}
	+\Vnom\left(\dyn\left(\state, \iter\right) - \step\pert\left(\state\right), \iter+1\right)
	 -  \Vnom\left(\dyn\left(\state, \iter\right), \iter +1\right)
	\notag\\
	&\stackrel{(b)}{\leq}  \!- \step \aa_3\!  \norm{\state}^2
	\!+\! \aa_4 \step  \!\norm{\pert\!\left(\state\right)}\left(\norm{\dyn\left(\state, \iter\right) \!-\! \step\pert\!\left(\state\right)} \!+\! \norm{\dyn\left(\state, \iter\right)}\right) 
	\notag\\
	&\stackrel{(c)}{\leq} - \step \aa_3\norm{\state}^2 + \step 2\aa_4\norm{\pert\left(\state\right)}\norm{\dyn\left(\state, \iter\right)} 
	+ \step^2 \aa_4\norm{\pert\left(\state\right)}^2
	\notag\\
	&\stackrel{(d)}{\leq} -  \step \aa_3 \norm{\state}^2 + \step 2\aa_4\aa_5 \epsilon \norm{\dyn\left(\state, \iter\right)} 
	+ \step^2 \aa_4\aa_5^2\epsilon^2,
	\label{eq:delta_V_original}
\end{align}
where in $(a)$ we add $\pm\Vnom\left(\dyn\!\left(\state, \iter\right), \iter + 1\right)$, in $(b)$ we use~\eqref{eq:delta_Vnom_negative} and~\eqref{eq:lipschitz_like_bound_Vnom} (cf. Lemma~\ref{lemma:nominal}), in $(c)$ we use the triangle inequality, while in $(d)$ we use~\eqref{eq:approx_bound}.
Now, since $\dyn(\state,\iter)$ is the sum of Lipschitz continuous functions (see~\eqref{eq:nominal_restricted_system} and Assumptions~\ref{ass:objective_functions},~\ref{ass:hat_lipschitz},~\ref{ass:NN_2}, and~\ref{ass:NN}), it is Lipschitz continuous too.
Further, we recall that $\dyn(0,\iter) = 0$ for all $\iter \in \N$.
Hence, by denoting with $\LF > 0$ its Lipschitz continuity constant, we bound~\eqref{eq:delta_V_original} as
\begin{align}
	\dVnom\left(\state, \iter\right) &\leq - \step \aa_3 \norm{\state}^2 + 2\step\aa_4 \aa_5  \LF \epsilon\norm{\state} + \step^2 \aa_4\aa_5^2\epsilon^2
	\notag\\
	&\stackrel{(a)}{\leq} - \step \aa_3 \norm{\state}^2 + \step \aa_4 \aa_5 \LF \pp \norm{\state}^2 
	+ \step \aa_4 \aa_5 \LF \epsilon^2/\pp + \step^2 \aa_4\aa_5^2 \epsilon^2
	\notag\\
	&\stackrel{(b)}{\leq}
	 -\step \aa_2\aat \norm{\state}^2  + \step  \aa_4 \aa_5( \LF/\pp + \step\aa_5) \epsilon^2
	\notag\\
	&\stackrel{(c)}{\leq} - \step \aat\Vnom\left(\state, \iter\right) + \step  \aa_4 \aa_5(\LF/\pp + \step\aa_5 ) \epsilon^2, \label{eq:delta_V_original_bound}
\end{align}
where in $(a)$ we apply the Young's inequality with parameter $\pp \in (0,(\aa_3 - \aa_2\aat)/(\aa_4\aa_5\LF))$ considering an arbitrarily chosen $\aat \in (0,\aa_3/\aa_2)$, in $(b)$ we rearrange the terms according to the choice of $\pp$, and in $(c)$ we use~\eqref{eq:positive_Vnom}.
Now, we arbitrarily fix $\step \in (0,\bar{\step})$,
define $\aa_6 \!\coloneq\! \frac{\aa_4 \aa_5}{\aat} (\frac{\LF}{\pp} +\step\aa_5)$ %
and bound~\eqref{eq:delta_V_original_bound} as
\begin{align}
	\dVnom\left(\state, \iter\right) &\leq - \step \aat \left(\Vnom\left(\state, \iter\right) - \epsilon^2\aa_6\right).\label{eq:inequality}
\end{align}
Because of the definition of $\dVnom$, the inequality~\eqref{eq:inequality} leads to
\begin{align}\label{eq:V_nom_update}
	\Vnom(\dyn(\state,\iter), \iterp) 
	& \leq \left(1 - \step \aat\right) \Vnom(\state, \iter) + \step \epsilon^2\aa_6 \aat.
\end{align}%
Thus, by iteratively applying~\eqref{eq:V_nom_update} 
and since $1 \!-\! \step\aat < 1$, we get
\begin{align}\label{eq:last}
	\Vnom\left(\statet, \iter\right) &\leq (1- \step \aat)^\iter\Vnom\left(\state\ud0, 0\right)  + \epsilon^2\aa_6.
\end{align}
The proof follows by combining~\eqref{eq:last} with~\eqref{eq:positive_Vnom} and $||\xt - \xstar|| \leq ||\statet||$ (see the definition of $\state$ in~\eqref{eq:state}), and by setting $\kk \coloneq \sqrt{\aa_2/\aa_1}$, $\kkk \coloneq \sqrt{1 - \step\aat}$, and $\bound \coloneq \sqrt{\aa_6/\aa_1}$.

\section{Numerical Simulations}
\label{sec:numerical_simulations}

\begin{figure}[htpb]
	\centering
	\subfloat[][]
	{\includegraphics[scale=\scaleplot]{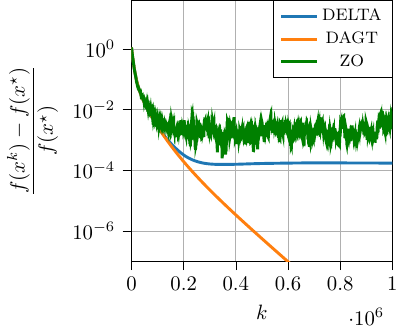} \ \ 
	\label{fig:cost}}
	\subfloat[][]
	{\includegraphics[scale=\scaleplot]{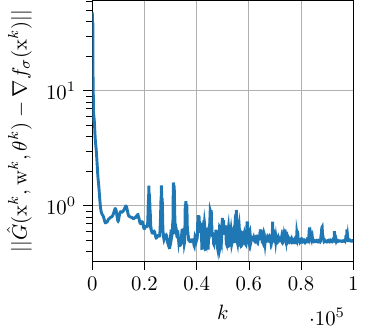}
	\label{fig:descent}
	}
	\caption{a) Comparison between \algo/, and DAGT~\cite{li2021distributed}, and DAGT with ZO gradient approximation by~\cite{flaxman2004online}. b) Evolution over time of $|| \hu(\xt,\wt,\tht) - \nabla \fs(\xt) ||$}
\end{figure}

In this section, we test \algo/ through numerical simulations.
We consider an instance of problem~\eqref{eq:aggregative_optimization_problem} with $N = 20$ agents.
In detail, for all $i \in \until{N}$, we consider the following setup.
We define the local cost functions $\f\du{i}$ as
\begin{align*}
	\f\du{i}(x\du{i},\sigma(x)) &\coloneq 
	\frac{1}{2}
	\begin{bmatrix}
		x\du{i}\\ \sigma(x)
	\end{bmatrix}\T
	\quadi
	\begin{bmatrix}
		x\du{i} \\ \sigma(x)
	\end{bmatrix}
	+ \lini\T
	\begin{bmatrix}
		x\du{i} \\ \sigma(x)
	\end{bmatrix}
	+ \expterm
	+ \constanti,
\end{align*}
where $x\du{i} \in \R$, $\sigma(x) \coloneq \sum_{i=1}^{N} \pi_i x\du{i} / N$.
Further, the parameters $\pi_i, \expa, \expb, \expc$ and the entries of $\quadi \in \R^{2\times 2}$ are drawn uniformly from $(0,1)$ enforcing $\quadi$ to be positive definite.
The entries of $\lini \in \R^{2}$ and $\constanti \in \R$ are drawn uniformly from $(0,20)$. 
Then, we set $\dxit(\iter)=\ampl\cos(2\pi\iter/\freq)$ and $\dsit(\iter)=\ampl\sin(2\pi\iter/\freq)$.
We pick local neural networks composed of $\layers$ layers with $\neurons$ neurons equipped with the softplus function.
As for the learning problem~\eqref{eq:learning_problem}, we consider the local cost functions
\begin{align*}
	\lossi(\uxi, \usi, \thi) \!\coloneq\! \frac{1}{2} \! (\f\du{i}(\uxi, \usi) \!-\! \hfi(\uxi, \usi, \thi))^2 \!\!\!+\! \norm{\thi}^2 \!\!\!.
\end{align*}
Further, we pick $x_i^0$ sampled from a Gaussian distribution, while $\thi^0$ is chosen by considering the Xavier uniform initialization~\cite{glorot2010understanding}.
We empirically tune $\step = \stepvalue$ and consider a random Erdős–Rényi graph with connectivity $p = 0.5$. 
Fig.~\ref{fig:cost} shows the evolution of the relative cost error $(\fs(\xt)-\fs(\xstar))/|\fs(\xstar)|$, where $\xstar$ is the optimizer of~\eqref{eq:aggregative_optimization_problem}.
We compare Algorithm~\ref{alg:delta} with an implementation of the Zeroth-Order (ZO) single-point method of~\cite{flaxman2004online} and with the Distributed Aggregative Gradient Tracking (DAGT) method of~\cite{li2021distributed}, which uses the exact gradients of $\f\du{i}$.
All simulations use the same stepsize and initial conditions.
All neural networks are implemented in TensorFlow, which provides built-in automatic differentiation~\cite{baydin2018automatic}.

\begin{figure}[htpb]
	\centering
	\includegraphics[scale=0.65]{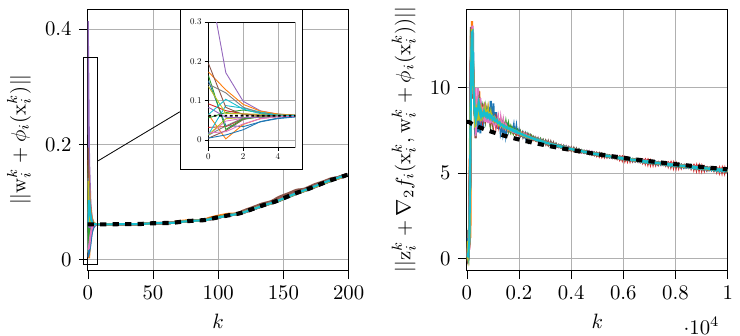}
	\caption{Evolution of the tracking dynamics. The dashed lines are the global $\sigma(\xt)$ (left) and $\sum_{j=1}^N \nabla_2 \hf\du{j}(\xjt,\sigma(\xt),\thjt)/N$ (right). The solid lines are their local agents estimates.}
	\label{fig:Tracker}
\end{figure}

\begin{figure}[htpb]
	\centering
	\vspace{-0.4cm}
	\subfloat[][Iteration $\iter = 0$.]
	   {\includegraphics[scale=\scalefig]{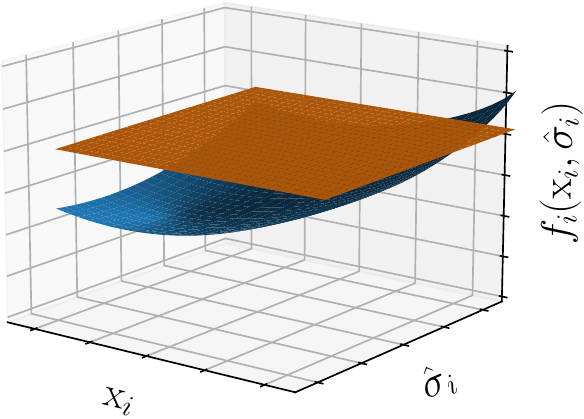}}\ \ \ \ 
	\subfloat[][Iteration $\iter = 20$.]
	   {\includegraphics[scale=\scalefig]{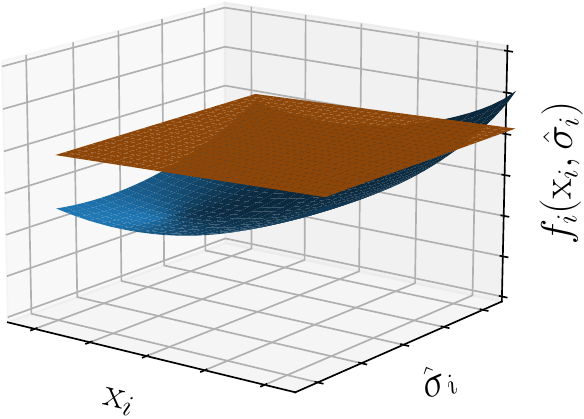}} \ \ \ \ 
	\subfloat[][Iteration $\iter = 500$.]
	   {\includegraphics[scale=\scalefig]{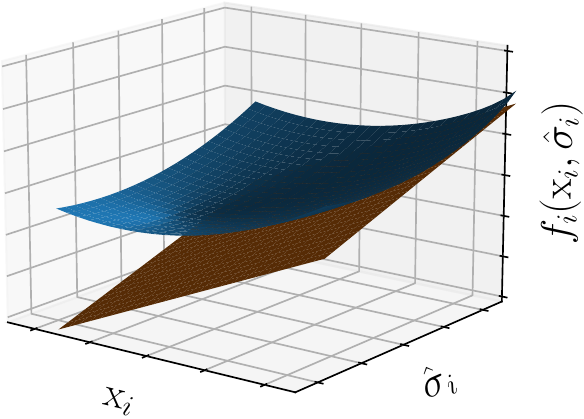}}\ \ \ \
	\subfloat[][Iteration $\iter = 10000$.]
	   {\includegraphics[scale=\scalefig]{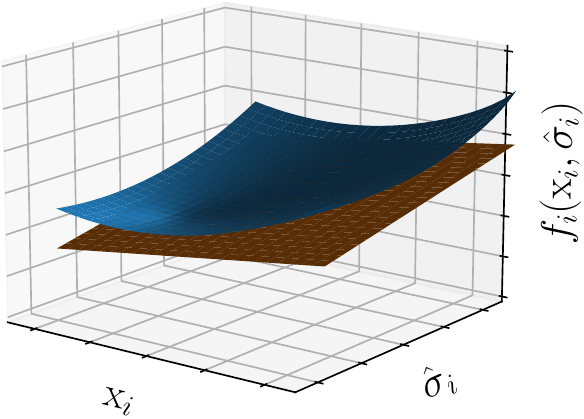}}
	   \caption{Agent $i$ perspective: evolution of the learned tangent space (in orange) over the iteration compared to $\f\du{i}$ (in blue).}
	   \label{fig:3d}
\end{figure}

\begin{figure}[H]
	\centering
	\subfloat[]
	{\includegraphics[scale=\scaleplot]{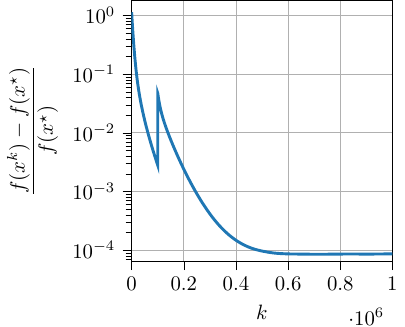}} \ \ \ \
	\subfloat[]
	{\includegraphics[scale=\scaleplot]{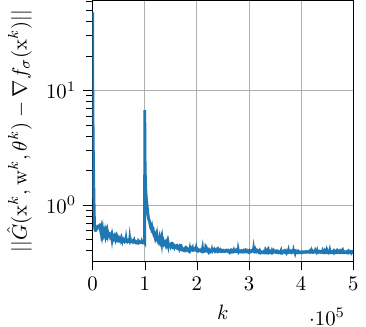}}
	\caption{Robustness of \algo/ to changes in the local costs. (a) Evolution of the relative cost error. (b) Evolution of the norm of the descent estimation error.}
	\label{fig:robustness}
\end{figure}

As predicted by Theorem~\ref{th:convergence}, Fig.~\ref{fig:cost} shows that the sequence $\{\xt\}_{\iter \in \N}$ generated by Algorithm~\ref{alg:delta} converges to a neighborhood of $\xstar$.
Then, Fig.~\ref{fig:descent} shows the evolution over time of the estimation error of the descent direction induced both by the neural networks and the tracking mechanism, i.e., $|| \hu(\xt,\wt,\tht) - \nabla \fs(\xt) ||$, where $\hu(\xt,\wt,\tht)$ is the estimated descent direction at iteration $\iter$ and thus it reads as
\begin{align*}
	\hu(\xt,\wt,\tht) &\coloneq \huxt(\xt, \wt + \phi(\xt), \tht) + \nabla \phi(\xt) (\hust\left(\xt, \wt + \phi(\xt), \tht\right) + \qt).
\end{align*} 
Fig.~\ref{fig:Tracker} illustrates the timescale separation in Algorithm~\ref{alg:delta}, 
showing that the consensus mechanism converges in fewer iterations than the optimization 
and learning components (cf. Figs.~\ref{fig:cost} and~\ref{fig:descent}).
In Fig.~\ref{fig:3d} we provide a graphical understanding of the learning capabilities of an agent. 
In detail, we show the evolution of the learned tangent space to the local cost $\f\du{i}$ around the current pair $(\xit,\wit + \phii(\xit))$.
Finally, we tested the robustness of \algo/ to changes in the local $\f_i$. In particular, at iteration $\iter = 10^5$ the parameters 
of each $\f\du{i}$ are perturbed by subtracting random values drawn from a uniform distribution in $(0,0.1)$. Figure~\ref{fig:robustness} shows that after this sudden change, the algorithm re-learns the new cost function gradients and converges to a neighborhood of the new optimal solution.

\section{Conclusions}

We proposed \algo/, a data-driven distributed algorithm for aggregative optimization in unknown scenarios.
\algo/ combines an optimization-oriented part with a tracking mechanism and a learning one.
The latter takes on local neural networks elaborating cost samples around the current solution estimates to approximate their gradients.
We analyzed \algo/ by using tools from system theory based on timescale separation and averaging theory.
In detail, we theoretically proved linear convergence in a neighborhood of the optimal solution whose size depends on the given accuracy capabilities of the neural networks.
Some numerical tests confirmed our findings.

\appendix

\setcounter{section}{0}
\renewcommand{\thesection}{\Alph{subsection}}

\subsection{Preliminaries on Averaging Theory}
\label{app:averaging}

We report~\cite[Theorem~2.2.2]{bai1988averaging}, which is a useful result in the context of averaging theory.
Consider the time-varying system
\begin{align} 
	\label{eq:discrete_time_system}
	\stateStp & = \chit + \step f(\chit, \step, \iter),
\end{align}
where $\chit \in \R^{\nc}$ is the system state, $f: \R^{\nc} \times \R \times \N \to \R^{\nc}$ describe its dynamics, and $\step > 0$ is a tunable parameter.
We now introduce a set of assumptions that allow for studying the stability properties of~\eqref{eq:discrete_time_system} by considering the so-called averaged system associated to~\eqref{eq:discrete_time_system}.
In particular, we impose these conditions within a sphere $\cB_r(0_{\nc})$ of radius $r > 0$.
\begin{assumption}\label{ass:f_av}
	The function $f$ is piecewise continuous in $\iter$ and the limit
	\begin{align}
		\label{eq:f_av}
		f\av(\chi) \coloneq \lim_{T \to \infty}\frac{1}{T} \sum_{\initer= \bar{\iter} +1}^{\bar{\iter} +T} f(\chi,0, \initer),
	\end{align}
	exists uniformly in $\bar{\iter} \in \N$ and $\chi \in \cB_r(0_{\nc})$.\oprocend
\end{assumption}
Once $f\av$ is well-posed, we can introduce the \emph{averaged system} associated to~\eqref{eq:discrete_time_system}, namely 
\begin{align}\label{eq:averaged_system}
	\chi\av\ud\iterp = \chi\av\ud\iter + \epsilon f\av(\chi\av\ud\iter),
\end{align}
with $\chit\av \in \R^{\nc}$.
Then, we enforce the following assumptions.
\begin{assumption}\label{ass:equilibrium}
	The origin is an equilibrium point for both~\eqref{eq:discrete_time_system} and~\eqref{eq:averaged_system}, namely 
	\begin{align*}
		f(0,\step,\iter) &= 0, \quad f\av(0) = 0,
	\end{align*}
	for all $\step \ge 0$ and $\iter \in \N$.\oprocend
\end{assumption}
\begin{assumption}\label{ass:lipschitz}
	There exist $l\du1, l\du2, l\av, \bar{\step}_1 > 0$ such that
	\begin{align*}
		\norm{f(\chi,\step,\iter) - f(\chi^\prime,\step,\iter)} &\leq l_1 \norm{\chi - \chi^\prime}\\
		\norm{f(\chi,\step,\iter) - f(\chi,\step^\prime,\iter)} &\leq l_2 \lvert \step - \step^\prime \rvert \norm{\chi}
		\\
		\norm{f\av(\chi) - f\av(\chi^\prime)} &\leq l\av \norm{\chi - \chi^\prime},
	\end{align*}
	for all $\step, \step^\prime \in \left(0 \right., \left. \bar{\step}_1\right]$, $\chi, \chi^\prime \in \cB_r(0_{\nc})$, and $\iter \in \N$. \oprocend 
\end{assumption}
\begin{assumption}\label{ass:nu}
	Let $\Delta f: \R^{\nc} \times \N \to \R^{\nc}$ be defined as
	\begin{align*}
		\Delta f(\chi, \iter) \coloneq f(\chi, 0, \iter) - f\av(\chi).
	\end{align*}
	Then, there exists a nonnegative strictly decreasing function $\nu(\iter)$ such that $\lim_{\iter \to \infty} \nu(\iter) = 0$ and
	\begin{subequations}%
		\begin{align*}
			\norm{\tfrac{1}{T}\textstyle\sum_{\initer = \bar{\iter} + 1}^{\bar{\iter} + T}\Delta f(\chi,\initer)} &\leq \nu(T)\norm{\chi}
			\\
			\norm{\tfrac{1}{T}\textstyle\sum_{\initer = \bar{\iter} + 1}^{\bar{\iter} + T}\frac{\partial\Delta f(\chi,\initer)}{\partial \chi}} &\leq \nu(T),
		\end{align*}
	\end{subequations}
	uniformly in $\bar{\iter} \in \N$ and $\chi \in \cB_{r}(0_{\nc})$.
	\oprocend
\end{assumption}
These assumptions ensure the following result.
\begin{theorem}\label{th:bai}\cite[Theorem~2.2.3]{bai1988averaging}
	Consider system~\eqref{eq:discrete_time_system} and let Assumptions~\ref{ass:f_av},~\ref{ass:equilibrium},~\ref{ass:lipschitz}, and~\ref{ass:nu} hold.
	If there exists $\bar{\step}_1 > 0$ such that, for all $\step \in (0,\bar{\step}_1)$, the origin is an exponentially stable equilibrium point of the averaged system~\eqref{eq:averaged_system}, then there exists $\bar{\step} \in (0,\bar{\step}_1)$ such that, for all $\step \in (0,\bar{\step})$, the origin is an exponentially stable equilibrium point of system~\eqref{eq:discrete_time_system}.\oprocend
\end{theorem}

\subsection{Proof of Lemma~\ref{lem:boundary_lyap}}
\label{sec:proof_bl}

	We point out that Assumption~\eqref{ass:graph} implies that $\Rmat\T\adjkron\Rmat$ is Schur. 
	Hence, for all $\Qmatfastentry  > 0$, there exist $\Plyapfast_1 = \Plyapfast_1\T > 0$ and $\Plyapfast_2 = \Plyapfast_2\T > 0$ solving the discrete-time Lyapunov equations
	\begin{subequations}\label{eq:Lyapunov_eq_fast}
		\begin{align}
			\label{eq:Lyapunov_eq_fast_1}
			\left(\Rmat\T \adjkron \Rmat\right)\T \Plyapfast_1 \left(\Rmat\T \adjkron \Rmat\right) - \Plyapfast_1 &= -\Qmatfastentry I_{(N-1)d}
			\\
			\label{eq:Lyapunov_eq_fast_2}
			\left(\Rmat\T \adjkron \Rmat\right)\T \Plyapfast_2 \left(\Rmat\T \adjkron \Rmat\right) - \Plyapfast_2 &=  - I_{(N-1)d},
		\end{align} 
	\end{subequations}
	where $\Qmatfastentry > 0$ will be fixed later.
	Let $\Plyapfast \coloneq \blkdiag(\Plyapfast_1, \Plyapfast_2) \in \R^{2(N-1)d}$ and $\Vfast: \R^{2(N-1)d}\to\R$ be defined as
	\begin{align*}
		\Vfast \left(\stateF\right) \coloneq
		\stateF\T
		\Plyapfast
		\stateF
		.
	\end{align*}
	Being $\Vfast$ quadratic, the conditions~\eqref{eq:positive_Vfast} and~\eqref{eq:lipschitz_like_bound_Vfast} are satisfied.
	To conclude the proof, we show that~\eqref{eq:delta_Vfast_negative} is verified studying the increment of $\Vfast$ along the trajectories of system~\eqref{eq:boundary_sys_compact}, i.e.,
	\begin{align}
		\dVfast (\stateF) &\coloneq \Vfast(\RARdiag \stateF + \gFast(\tx, \stateF, \tth)) - \Vfast(\stateF)
		\notag\\
		&= 
		\stateF\T
		(\RARdiag\T \Plyapfast \RARdiag - \Plyapfast)
		\stateF
		+ 2 \stateF\T \RARdiag\T \Plyapfast \gFast(\tx, \stateF, \tth)
		+
		\Delta\hust(\xt, \stateF_1, \tht)\T
		\Kmat
		\Delta\hust(\xt, \stateF_1, \tht),\label{eq:increment_1}
	\end{align}
	where 
	$\Kmat \coloneq \begin{bmatrix}
		0& \Rmat\T\adjminusI
	\end{bmatrix}\Plyapfast
		\begin{bmatrix}
			0&
			\Rmat\T\adjminusI
		\end{bmatrix}\T$. 
	By looking at the definition of $\hust$ (cf.~\eqref{eq:ust}) and $\Delta\hust$ (cf.~\eqref{eq::u2}), since $\nabla_2 \hfi$ is Lipschitz continuous (cf. Assumption~\ref{ass:hat_lipschitz}), we have
	\begin{align}
		\|\Delta\hust(\tx \!+\! \xstar\!, \twort, \tth \!+\! \thstar\!)\| \leq \Lfsigma \|\Rmat \twort\| ,\label{eq:lip_u2}
	\end{align}
	for all $\twort \in \R^{(N-1)d}$.
	Now, we use~\eqref{eq:Lyapunov_eq_fast} to bound~\eqref{eq:increment_1} as
	\begin{align}
		\dVfast\left(\stateF\right) 
		&=
		- \Qmatfastentry \norm{\stateF_1}\ud{2} -  \norm{\stateF_2}\ud{2} + 2 \stateF\T \RARdiag\T \Plyapfast \gFast(\tx, \stateF, \tth) 
		\notag\\
		&\hspace{12pt} 
		+ \Delta\hust(\x, \stateF_1, \tth)\!\T\!(\Rmat\T\adjminusI)\!\T\!\Plyapfast_2 (\Rmat\T\adjminusI){\Delta}\hust(\x, \stateF_1, \tth)
		\notag\\
		&\stackrel{(a)}{\leq}
		- \! \Qmatfastentry \norm{\stateF_1}\ud{2} \!-\!  \norm{\stateF_2}\ud{2} 
		\! + \! 2 \norm{\stateF\du{2}} \norm{\Rmat\T \adjkron \Rmat \Plyapfast\du{2} \Rmat \adjkron}\!\norm{\Delta\hust(\x, \stateF_1, \tth)}
		\notag\\
		&\hspace{12pt} + \norm{(\Rmat\T\adjminusI)\!\T\!\Plyapfast_2 (\Rmat\T\adjminusI)} \norm{\Delta\hust(\x, \stateF_1, \tth)}\ud{2}
		\notag\\
		&\stackrel{(b)}{\leq} 
		\!- \Qmatfastentry \norm{\stateF_1}\ud{2} -  \norm{\stateF_2}\ud{2}
		+2 \Lfsigma \norm{\Rmat\T \adjkron \Rmat \Plyapfast\du{2} \Rmat \adjkron}  \norm{\Rmat \stateF_1} \norm{\stateF\du{2}}
		\notag\\
		&\hspace{12pt}
		+ \norm{(\Rmat\T\adjminusI)\!\T\!\Plyapfast_2 (\Rmat\T\adjminusI)}\Lfsigma\ud{2} \norm{\Rmat \stateF_1}\ud{2}
		,
		\label{eq:increment_2}
	\end{align}
	where in $(a)$ we use the Cauchy-Schwarz inequality and the definition of $\gFast$ (cf.~\eqref{eq:gFast}), in $(b)$ we use equation~\eqref{eq:lip_u2}. %
	Let us define $\kappa_1 \coloneqq \norm{(\Rmat\T\adjminusI)\!\T\!\Plyapfast_2 (\Rmat\T\adjminusI)}\Lfsigma\ud{2} \norm{\Rmat}\ud{2}$ and $\kappa_{12} \coloneq \Lfsigma \norm{\Rmat\T \adjkron \Rmat \Plyapfast\du{2} \Rmat \adjkron}  \norm{\Rmat}$.
	Then, we compactly rewrite~\eqref{eq:increment_2} as
	\begin{align}
		\label{eq:increment_3}
		\dVfast\left(\stateF\right) 
		&\leq 
		- \begin{bmatrix}
			\norm{\stateF_1}
			\\ 
			\norm{\stateF_2}
		\end{bmatrix}\T
		\underbrace{\begin{bmatrix}
			\Qmatfastentry - \kappa_1  & \kappa_{12}
			\\
			\kappa_{12}  &  1
		\end{bmatrix}}_{=:\Pmatfast}
		\begin{bmatrix}
			\norm{\stateF_1}
			\\
			\norm{\stateF_2}
		\end{bmatrix}.
	\end{align}
	By Sylvester Criterion, we know that $\Pmatfast > 0$ if and only if 
	\begin{align}\label{eq:conditions}
		\Qmatfastentry > \kappa_1 \qquad \text{and} \qquad \Qmatfastentry > \kappa_1  + \kappa_{12}^2.
	\end{align}	
	We set $\Qmatfastentry > \bar{\Qmatfastentry}_1 \coloneq  \kappa_1 +\kappa_{12}^2$ and, thus, we bound~\eqref{eq:increment_3} as
	\begin{align*}
		\dVfast\left(\stateF\right)&\leq - \lambda\du{\text{min}}\left(\Pmatfast\right) \norm{\stateF}\ud{2},
	\end{align*}
	where $\lambda\du{\text{min}}\left(\Pmatfast\right)>0$ denotes the smallest eigenvalue of $\Pmatfast$.
	Hence, also~\eqref{eq:delta_Vfast_negative} is achieved and, thus, the proof concludes.

\subsection{Proof of Lemma~\ref{lemma:av}}
\label{sec:proof_av}

By relying on the accuracy characterization~\eqref{eq:theta_star_nabla_ell}, we get
\begin{align}
\nablalossav\left(\tx, 0\right) = \tfrac{1}{\Tp}\textstyle\sum_{\iter=1}^{\Tp} \pertz\left(\tx + \xstar,\iter\right).\label{eq:sum_grad_d}
\end{align}
for all $\tx \in \R^{\n}$.
Since $\thstar$ is a minimizer of problem~\eqref{eq:learning_problem} for $\iter = \Tp$ (cf. Assumption~\ref{ass:NN_2}), it holds $\nablalossav\left(\tx, 0\right) = 0$, which, combined with~\eqref{eq:sum_grad_d}, leads to
	$\sum_{\iter=1}^{\Tp} \pertz\left(\tx + \xstar,\iter\right)= 0$
for all $\tx \in \R^{\n}$.
Thus, system~\eqref{eq:reduced_sys_av_with_p} reduces to 
\begin{subequations}
	\label{eq:reduced_sys_av}
	\begin{align}
		\tthtpav &= \tthtav - \step \nablalossav(\txav,\tthtav)
		\\
		\txtpav &= \txtav - \step \hud(\txtav,\tthtav).
	\end{align}
\end{subequations}
Now, let $\Vslow: \R^{\dimtheta} \times \R^{\n} \to \R$ be the Lyapunov function
\begin{align*}
	\Vslow(\tthav,\txav) = \ath \norm{\tthav}^2 + \norm{\txav}^2,
\end{align*}
where $\ath > 0$ will be fixed in the following.
Let us introduce the increments $\Delta\Vslowth(\tthav, \txav)$ and $\Delta\Vslowtx(\tthav, \txav)$ defined as
\begin{align*}
	\Delta\Vslowth(\tthav, \txav) &\coloneq \|\tthav - \step \nablalossav(\txav, \tthav)\|^2 - \|\tthav\|^2
	\\
	\Delta\Vslowtx(\tthav, \txav) &\coloneq \|\txav - \step \hud (\txav, \tthav)\|^2 - \|\txav\|^2.
\end{align*}
Then, by considering $(\tthav, \txav) \in \R^{\dimtheta} \times \R^{\n}$ with $(\tthav + \thstar) \in \thSet$, the increment of $\Vslow$ along the trajectories of system~\eqref{eq:reduced_sys_av} is
\begin{equation}
	\label{eq:delta_Vslow}
	\begin{aligned}
		\dVslow(\tthav, \txav) %
		 =  \ath \Delta\Vslowth(\tthav, \txav)  +  \Delta\Vslowtx(\tthav, \txav).
	\end{aligned}
\end{equation}
As for $\Delta\Vslowth(\tthav, \txav)$, we expand the square norm and obtain
\begin{align}
	\Delta\Vslowth(\tthav, \txav)  
	&=  -  \! 2\step\nablalossav(\txav, \tthav)\T\tthav 
	\!+\! \step^2\norm{\nablalossav(\txav, \tthav)}^2
	\notag\\ 
	&\stackrel{(a)}{\leq} 
	-\step \tfrac{2\multh \Llth}{\multh \!+\! \Llth}\norm{\tthav}^2 \!-\! \step  (\tfrac{2}{\multh \!+\! \Llth} \!-\! \gamma)\norm{\nablalossav(\txav, \tthav)}^2,
	\label{eq:Vth_bound}
\end{align}
where in $(a)$ we use $\nablalossav(\txav, 0)=0$, the $\multh$-strong convexity of $\lossav$ in $\thSet$, and the $\Llth$-Lipschitz continuity of $\nablalossav$ (cf. Assumption~\ref{ass:NN}).
Then, we expand $\Delta\Vslowtx$ and obtain
\begin{align}
	\Delta\Vslowtx(\tthav, \txav) 
	\!&=\!
	-2 \step \hud(\txav, \tthav)\T \txav + \step^2 \norm{\hud(\txav, \tthav)}^2
	\notag\\
	&\!\stackrel{(a)}{=}  - 2 \step \nabla\fs(\txav + \xstar)\T \txav
	\notag\\
	& \hspace{12pt} -2 \step (\hud(\txav, \tthav) \!-\! \hud(\txav, 0))\T\! \txav
	\notag\\
	& \hspace{12pt} + \step^2 \!\norm{\hud(\txav,\! \tthav) \!-\! \hud(\txav,\! 0) \!+\! \hud(\txav,\! 0)}^2
	\notag\\
	&\!\stackrel{(b)}{=}
	-2 \step \nabla\fs(\txav + \xstar)\T \txav + 2 \step^2 \norm{\hud(\txav, 0)}^2\!
	\notag\\
	& \hspace{13pt} 
	\!+ 2 \step \!\norm{\hud(\txav, \tthav) \!-\! \hud(\txav, 0)} \norm{\txav}
	\notag\\
	& \hspace{13pt} \!+\! 2 \step^2 \!\norm{\hud(\txav, \tthav) \!-\! \hud(\txav, 0)}^2
	,\label{eq:Vx_bound}
\end{align}
where in $(a)$ we add $\pm\nabla\fs(\txav+\xstar) \!=\! \hud(\txav, 0)$, and $(b)$ uses the triangle and Cauchy-Schwarz inequalities and $\nabla\fs(\xstar) \!=\! \hud(0,0) \!=\! 0$. 
By Lipschitz continuity of $\nabla_1 \hfi$, $\nabla_2 \hfi$, and $\phii$ (cf. Assumptions~\ref{ass:objective_functions} and~\ref{ass:hat_lipschitz}), there exists $\Lfth \!>\! 0$ such that 
\begin{align}\label{eq:bound_lip_hud}
	\|\hud(\txav, \tthav) - \hud(\txav, \bar{\thh}\av)\| \leq \Lfth\|\tthav - \bar{\thh}\av\|,
\end{align}
for all $\txav \!\in\! \R^{\n}$ and $\tthav, \bar{\thh}\av \!\in\! \R^{\dimtheta}$.
By~\eqref{eq:bound_lip_hud}, $\mufx$-strong convexity of $\fs$, and $\Lfx$-Lipschitz continuity of $\nabla \fs$, we bound~\eqref{eq:Vx_bound} as
\begin{align}
	\Delta\Vslowtx(\tthav, \txav) %
	&\leq - \step \tfrac{2 \mufx \Lfx}{\mufx + \Lfx} \norm{\txav}^2 - \step (\tfrac{2}{\mufx + \Lfx} - \step) \|\hud(\txav, 0)\|^2
	\notag\\
	&\hspace{11pt}+ 2 \step \Lfth \norm{\txav} \|\tthav\| + 2 \step^2 \Lfth^2 \|\tthav\|^2.
	\label{eq:Vx_final_bound}
\end{align} 
Let us arbitrarily set $\betax \!\in\! (0, \frac{2 \mufx \Lfx}{\mufx + \Lfx})$ and $\betath \in (0, \frac{2\multh \Llth}{\multh + \Llth})$ and define $\bar{\step}_2 \coloneq  \min \left\{\frac{\betax}{\mufx \Lfx}, \frac{\betath}{\multh \Llth}\right\}$. 
Then, by choosing $\step \in (0,\bar{\step}_2)$ and using~\eqref{eq:Vth_bound}-\eqref{eq:Vx_final_bound}, we bound $\dVslow$ (cf.~\eqref{eq:delta_Vslow}) as
\begin{align}
	\dVslow(\tthav, \txav) \! \leq \! - \!
	\begin{bmatrix}
		\norm{\tthav} \\ \norm{\txav} 
	\end{bmatrix}\T\!\!\!
	\underbrace{\begin{bmatrix}
		\ath \step \betath - 2 \step^2 \Lfth^2  & -\step \Lfth\\
		-\step \Lfth & \step \betax  
	\end{bmatrix}}_{=:\Pmatslow(\ath)}
	\!\!
	\begin{bmatrix}
		\norm{\tthav} \\ \norm{\txav} 
	\end{bmatrix}\!.\label{eq:dVslow}
	\notag
\end{align}
Choose $\ath \!>\! \bar{\ath} \!\coloneq\!(\Lfth^2 \left(1 \!+\! 2\step \betax\right))/\betath\betax$.
By Sylvester Criterion, we get $\Pmatslow(\ath) > 0$ and, thus, 
we bound $\dVslow(\stateSav)$ as
\begin{align*}
	\dVslow\left(\stateSav\right)\leq -\lambda\du{\text{min}}\left(\Pmatslow(\ath)\right) \norm{\stateSav}\ud{2}	
\end{align*}
where $\lambda\du{\text{min}}(\Pmatslow(\ath)) > 0$ represents the smallest eigenvalue of $\Pmatslow$ and $\stateSav \coloneq \col(\tthav, \txav)$. 
Thus, the exponential stability of the origin is proved (see, e.g.~\cite[Th.~13.11]{haddad2008nonlinear}).

\subsection{Proof of Lemma~\ref{lem:reduced_stability}}
\label{sec:proof_reduced_stability}

The proof relies on Theorem~\ref{th:bai} (cf. Appendix~\ref{app:averaging}).
Hence, we need to check that Assumptions~\ref{ass:f_av},~\ref{ass:equilibrium},~\ref{ass:lipschitz}, and~\ref{ass:nu} are satisfied. %
First, Assumption~\ref{ass:f_av} is guaranteed by the periodicity of $\dxit$ and $\dsit$ (cf. Assumption~\ref{ass:NN_2}).
Second, Assumption~\ref{ass:equilibrium} (i.e., the fact that the origin is an equilibrium point of both~\eqref{eq:reduced_sys_imp} and the averaged system~\eqref{eq:reduced_sys_av}) is satisfied in light of~\eqref{eq:theta_star_nabla_ell} in Assumption~\ref{ass:NN} and since $\thstar$ is the solution to~\eqref{eq:learning_problem} with $\iter = \Tp$ (cf. Assumption~\ref{ass:NN_2}).
Third, the Lipschitz properties in Assumption~\ref{ass:lipschitz} follow by Assumptions~\ref{ass:NN_2} and~\ref{ass:NN}. 
Fourth, the periodicity of $\dxit$ and $\dsit$ (cf. Assumption~\ref{ass:NN}) allows for satisfying Assumption~\ref{ass:nu} with $\nu\left(T\right) \coloneq \frac{L}{1+T}$ where $L>0$ depends on the Lipschitz constants of $\nabla_3 \loss$ and $\nabla_3 \lossav$.
Finally, Lemma~\ref{lemma:av} proves that the origin is an exponentially stable equilibrium of the averaged system~\eqref{eq:reduced_sys_av}.
Hence, by Theorem~\ref{th:bai}, there exist $\bar{\step}_1 \in (0,\bar{\step}_2)$ and $r > 0$ such that, for all $\step \in (0,\bar{\step}_1)$ and $\stateS\ud0\in \cX$, the origin is an exponentially stable equilibrium of~\eqref{eq:reduced_sys_imp}.
The proof follows by the Converse Lyapunov Theorem (see, e.g.,~\citeConverse).

\subsection{Proof of Lemma~\ref{lemma:nominal}}
\label{sec:proof_nominal}

The proof is based on timescale separation.
First, Lemma~\ref{lem:boundary_lyap} provides a Lyapunov function proving the global exponential stability of the origin of~\eqref{eq:boundary_sys_compact}, i.e., the boundary-layer system associated to~\eqref{eq:nominal_restricted_system}.
Second, Lemma~\ref{lem:reduced_stability} provides a Lyapunov function proving exponential stability of the origin of~\eqref{eq:reduced_sys_compact}, i.e., the reduced system associated to~\eqref{eq:nominal_restricted_system}.
Third, the dynamics of system~\eqref{eq:nominal_restricted_system} and the equilibrium map~\eqref{eq:parametrized_equilibria} are Lipschitz continuous (see Assumptions~\ref{ass:objective_functions},~\ref{ass:hat_lipschitz},~\ref{ass:NN_2}, and~\ref{ass:NN}).
We thus guarantee exponential stability of $(\thstar, \xstar, \h(\thstar, \xstar))$ for system~\eqref{eq:nominal_restricted_system} by slightly extending~\cite[Th.~II.3]{carnevale2023admm}.
In detail, the mentioned result provides global results but requires global results in both the boundary-layer and reduced system, while Lemma~\ref{lem:reduced_stability} (i.e., the stability result of the reduced system) holds only locally.
Hence, we conclude that the origin is an exponentially stable equilibrium point of system~\eqref{eq:nominal_restricted_system} for any initial condition  $(\tth\ud0, \tx\ud0, \wort\ud0,\qort\ud0) \in \domVpert$.
The proof follows by the Converse Lyapunov Theorem, see, e.g.,~\citeConverse.

\bibliographystyle{ieeetr}
\bibliography{bibliography}

\end{document}